\documentclass[12pt]{article}

\usepackage{yhmath}
\usepackage{amstext    }
\usepackage{amsthm    }
\usepackage{a4}
\usepackage[mathscr]{eucal}
\usepackage{mathrsfs}

\usepackage{amsmath}
\usepackage{amssymb}
\usepackage{amscd}

\newtheorem{theorem}{Theorem}[section]
\newtheorem{definition}[theorem]{Definition}
\newtheorem{proposition}[theorem]{Proposition}
\newtheorem{corollary}[theorem]{Corollary}
\newtheorem{lemma}[theorem]{Lemma}

\newcommand{\cali}[1]{\mathscr{#1}}

\newcommand{\lof}{\mathop{\mathrm{{log^\star}}}\nolimits}
\newcommand{\dist}{\mathop{\mathrm{dist}}\nolimits}
\newcommand{\diam}{\mathop{\mathrm{diam}}\nolimits}

\newcommand{\id}{{\rm id}}

\newcommand{\reg}{{\rm reg}}
\newcommand{\sing}{{\rm sing}}

\newcommand{\radius}{\mathop{\mathrm{radius}}\nolimits}

\newcommand{\Cc}{\cali{C}}

\newcommand{\Lc}{\cali{L}}

\newcommand{\Vc}{\cali{V}}

\newcommand{\C}{\mathbb{C}}
\newcommand{\D}{\mathbb{D}}
\renewcommand{\H}{\mathbb{H}}

\newcommand{\Z}{\mathbb{Z}}
\newcommand{\R}{\mathbb{R}}
\newcommand{\T}{\mathbb{T}}

\newcommand{\U}{\mathbb{U}}

\renewcommand{\P}{\mathbb{P}}

\title{Entropy for hyperbolic Riemann surface laminations II}

\author{Tien-Cuong Dinh, Viet-Anh Nguy{\^e}n and Nessim Sibony}

\begin{document}

\maketitle

\begin{abstract}
Let $(X,\Lc,E)$ be a Brody hyperbolic foliation by Riemann surfaces with linearizable isolated singularities on a compact complex surface. We show that its hyperbolic entropy is finite. We also estimate the modulus of continuity of the Poincar\'e metric on leaves. The estimate holds for foliations on manifolds of higher dimension. 
\end{abstract}

\noindent
{\bf Classification AMS 2010:} 37F75, 37A.


\section{Introduction} \label{introduction}

In this second part, we study  Riemann surface foliations with tame singular points.  
We say that a holomorphic vector field $F$ in $\C^k$ is {\it generic linear} if
it can be written as 
$$F(z)=\sum_{j=1}^k \lambda_j z_j {\partial\over \partial z_j}$$
where the $\lambda_j$ are non-zero complex numbers.
The integral curves of $F$ define a Riemann surface foliation in
$\C^k$. 
The condition $\lambda_j\not=0$ for every $j$  implies that the foliation   
has an isolated singularity at 0. 

Consider a Riemann surface foliation with singularities $(X,\Lc,E)$  in a  complex manifold $X$. We assume that
the singular set $E$ is discrete. By {\it foliation}, we mean that $\Lc$ is transversally holomorphic.
We say that a singular point $e\in E$  is {\it linearizable} if  there are  
local holomorphic coordinates on a neighborhood $\U$ of $e$ in which 
the foliation  is  given by  a generic linear vector field.
 
The purpose  of   the  present paper is to study the notion of  Brody  hyperbolicity  for compact Riemann surface foliations that will be given in 
Definition  \ref{uniform_hyperbolic_laminations} below.  
Here is our main result.

\begin{theorem} \label{th_main}  Let $(X,\Lc,E)$ be a singular foliation by
Riemann surfaces  on a compact complex  surface $X$. Assume that the singularities are linearizable and that
 the foliation is  Brody hyperbolic.
Then, its hyperbolic entropy $h(\Lc)$ is  finite.
\end{theorem}

We deduce from the above theorem and a result by Glutsyuk \cite{Glutsyuk} and Lins Neto \cite{Neto} the  following corollary. It can be applied to foliations of degree at least 2 with hyperbolic singularities.

\begin{corollary} \label{cor_main}  
Let $(\P^2,\Lc,E)$ be a singular foliation by
Riemann surfaces  on the complex projective plane $\P^2$. Assume that the singularities are linearizable. 
Then, the hyperbolic entropy $h(\Lc)$ of $(\P^2,\Lc,E)$ is  finite.  
\end{corollary}

In Section \ref{section_local}, we prove the finiteness of the entropy in the local setting near a singular point in any dimension. The main result in this section does not imply Theorem \ref{th_main}. However, it clarifies a difficulty due to singular points. We use here a division of a neighbourhood of a singular point into adapted cells. The construction of these cells is crucial in the proof of Theorem \ref{th_main}.

In Section \ref{section_poincare}, we will estimate the modulus of continuity for the Poincar\'e metric along the leaves of the foliation. We will use the notion of conformally $(R,\delta)$-close maps as in the case without singularities but there are several technical problems. Indeed, we have to control the phenomenon that leaves may go in and out singular flow boxes without any obvious rule. 

The estimates we get in Section \ref{section_poincare}  hold uniformly on $X$. They are still far from being sufficient in order to get Theorem \ref{th_main}, i.e.
the finiteness of entropy in the global setting. 
The proof of that result is more delicate. It is developed in the last three sections. 
We are only able to handle the interaction of the two difficulties mentioned above for foliations on complex surfaces and we conjecture that our main result is true in any dimension. A basic idea is the use of conformally $(R,\delta)$-close maps from leaves to leaves with small Beltrami's coefficient as in the case without singularities. 
In particular, we will glue together local orthogonal projections from leaves to leaves. 
However, we have to face the problem that the  Poincar\'e metric is not bounded from above by a smooth Hermitian metric on $X$ and hence the Beltrami coefficient of orthogonal projections from leaves to leaves is not small near the singularities. We solve this difficulty by replacing these projections near the singularities by adapted 
holomorphic maps which exist thanks to the nature of the singular points.

We can introduce, for any extremal harmonic measure $m$, the metric entropy $h(m)$,
the local entropies $h^\pm$ and the transverse local entropies $\widetilde h^\pm$. As in the case without singularities, we can show that $\widetilde h^\pm$ are constant $m$-almost everywhere and $h^\pm\leq \widetilde h^\pm+2$. We believe that $h^\pm=\widetilde h^\pm+2$ and then $h^\pm$ are constant $m$-almost everywhere but the question is still open. Other open problems stated for laminations without singularities in \cite{DinhNguyenSibony2} can be also considered for singular foliations.

The rest of the paper is quite technical. 
In order to simplify the presentation, we will not try to get sharp constants in our estimates and we prefer to give simple statements. In particular, for the results below, we assume that the hyperbolic time $R$ is large enough and we heavily use the fact that $R$ is larger than any fixed constants. Our notation is carefully chosen and we invite the reader to keep in mind the following remarks and conventions.

\bigskip
\noindent
{\bf Main notation.} We use the same notation as in Part 1 of this paper \cite{DinhNguyenSibony2}, e.g.   
$\D$, $r\D$, 
$\D_R$, $\omega_P$, $\dist_P$, $\diam_P$, $\eta$, $L_x$, $\phi_x:\D\to L_x$, $\dist_R$, $\dist_{\D_R}$ and $\mu_\tau$. 
Denote by $\D(\xi, R)$ the disc of center $\xi$ and of radius $R$ in the Poincar\'e disc $\D$ and $\lof(\cdot):=1+|\log(\cdot)|$ a log-type function.
Conformally $(R,\delta)$-close points are defined as in \cite{DinhNguyenSibony2} with an adapted constant $A:=3c_1^2$ and for $\delta\leq e^{-2R}$. 
If $D$ is a disc, a polydisc or a ball of center $a$ and $\rho$ is a positive number, then $\rho D$ is the image of $D$ by the homothety $x\mapsto \rho x$, where $x$ is an affine coordinates system centered at $a$. In particular,
$\rho\U_i$ and $\rho \U_e$ correspond respectively to a regular flow box and to a flow box at a singular point $e\in E$.

\bigskip
\noindent
{\bf Flow boxes and metric.} We only consider flow boxes which are biholomorphic to $\D^k$. For regular flow boxes, i.e. flow boxes outside the singularities, the plaques are identified with the discs parallel to the first coordinate axis. Singular flow boxes are identified to their models described in Section \ref{section_local}. 
In particular, the leaves in a singular flow box are parametrized in a canonical way using holomorphic maps $\varphi_x:\Pi_x\to L_x$, where $\Pi_x$ is a convex polygon in $\C$. 

For each singular point $e\in E$, we fix a singular flow box $\U_e$ such that $3\U_e\cap 3\U_{e'}=\varnothing$ if $e\not=e'$. We also cover $X\setminus \cup {1\over 2}\U_e$ by regular flow boxes ${1\over 2}\U_i$ which are fine enough. In particular, each $\U_i$ is contained in a larger regular flow box $2\U_i$ with $2\U_i\cap {1\over 4}\U_e=\varnothing$ such that if $2\U_i$ intersects $\U_e$, it is contained in $2\U_e$. More consequences of the small size of $\U_i$ will be given when needed. We identify $\{0\}\times \D^{k-1}$ with a transversal $\T_i$ of $\U_i$ and call it {\it the distinguished transversal}. For $x\in\U_i$, denote by $\Delta(x,\rho)$ the disc of center $x$ and of radius $\rho$ contained in the plaque of $x$.

Fix a Hermitian metric $\omega$ on $X$ which coincides with the standard Euclidean metric on each singular flow box $2\U_e\simeq 2\D^k$.

\bigskip
\noindent
{\bf Other notation.}
Denote by $L_x[\epsilon]$ the intersection of $L_x$ with the ball of center $x$ and of radius $\epsilon$ with respect to the metric induced on $L_x$ by the Hermitian metric on $X$. 
It should be distinguished from $L_x(\epsilon):=\phi_x(\D_\epsilon)$. 
We only use $\epsilon$ small enough so that $L_x[\epsilon]$ is connected and simply connected. Fix a constant $\epsilon_0>0$ small enough so that if $x$ is outside the singular flow boxes ${1\over 2}\U_e$, then $L_x[\epsilon_0]$ is contained in a plaque of a regular flow box ${1\over 2}\U_i$. If
$L_x(\epsilon_0):=\phi_x(\D_{\epsilon_0})$ is not contained in a singular flow box ${1\over 4}\U_e$, then $\phi_x$ is injective on $\D_{\epsilon_0}$.

The constant $\gamma$ is given in Lemma \ref{lemma_J_l_close}. The constant $\lambda$ in Section \ref{section_local} can be equal to $\lambda_*$ but we will use later the case with a large constant $\lambda$. Define 
$\alpha_1:=e^{-e^{7\lambda R}}$ and $\alpha_2:=e^{-e^{23\lambda R}}$.  
The set $\Sigma$ and the maps $J_l$, $\Phi_{x,y}$, $\Psi_{x,y}$,$\widetilde\Psi_{x,y}$ are introduced in Section \ref{section_local}.
The constants $A$ and $c_i$ are introduced in Section \ref{section_poincare} with $A:=3c_1^2$. The constants $m_0,m_1,\hbar,p$ are fixed just before Lemma \ref{lemma_p_points_1}, the constant $t$ just after that lemma. They satisfy $c_1\ll m_0$, $c_1m_0\ll m_1$, $m_1\hbar\ll\epsilon_0$ and $p:=m_1^4$.

\bigskip
\noindent
{\bf Orthogonal projections.} For $x=(x_1,\ldots, x_k)\in \C^k$, define the norm $\|x\|_1$ as
$\max |x_j|$. 
We choose $\epsilon_0>0$ small enough so that if $x$ and $y$ are two points outside the singular flow boxes ${1\over 4}\U_e$ such that $\dist(x,y)\leq \epsilon_0$, then for the Euclidean metric, the local orthogonal projection $\Phi$ from $L_x[\epsilon_0]$ to $L_y[10\epsilon_0]$ is well-defined and its image is contained in $L_y[3\epsilon_0]$. In a singular flow box $\U_e\simeq \D^k$, since the metric and the foliation are invariant under homotheties, when $\dist(x,y)\leq \epsilon_0\|x\|_1$, the local orthogonal projection is well-defined from $L_x[\epsilon_0\|x\|_1]$ to $L_y[10\epsilon_0\|x\|_1]$ with image in 
$L_y[3\epsilon_0\|x\|_1]$. 
If moreover, $x,y$ are very close to each other and are outside the coordinate hyperplanes, 
a global orthogonal projection $\Phi_{x,y}$ from $L_x$ to $L_y$ is constructed in Lemma \ref{lemma_proj_E}. 
All projections described above are called {\it the basic projections} associated to $x$ and $y$. They are not holomorphic in general. 

In order to construct a map $\psi$ satisfying the definition of conformally $(R,\delta)$-close points, we have to glue together basic projections. In the case without singularities \cite{DinhNguyenSibony2}, we have carefully shown that the gluing is possible, i.e.
there is no monodromy problem. The same arguments work in the case with singularities. We sometimes skip the details on this point in order to simplify the presentation.


\section{Local models for singular points} \label{section_local}

In this section, we give a description of the local model for linearizable singularities. We also prove the finiteness of entropy in this setting. The construction of cells and other auxiliary results given at the end of the section will be used in the proof of Theorem \ref{th_main}.

Consider the foliation $(\D^k,\Lc,\{0\})$ which is the restriction to $\D^k$ of the foliation associated to the vector field
$$F(z)=\sum_{j=1}^k \lambda_jz_j{\partial\over \partial z_j}$$
with $\lambda_j\in\C^*$.  The foliation is singular at the origin. We use here the Euclidean metric on $\D^k$.
The notation $\widehat L_x$, $\widehat\phi_x$, $\widehat\eta$, $\dist_R$ and $h(\cdot)$  below are defined as in the case of general foliations. Here, we use a hat for some notations in order to avoid the confusion with the analogous notations that we will use later in the global setting.

Define
$$\lambda_*:={\max \left\lbrace |\lambda_1|,\ldots,|\lambda_k| \right\rbrace \over    \min \left\lbrace |\lambda_1|,\ldots,|\lambda_k| \right\rbrace}\cdot$$

\begin{theorem}\label{th_local_entropy}
For every compact subset $K$ of  $\D^k$, the  hyperbolic  entropy $h(K)$ of $K$ is bounded from above by
$70k \lambda_*$.
\end{theorem}

Observe that the entropy of $K$ is bounded independently of $K$.
For the proof of this result, we will construct a  division of $\D^k$ into cells whose shape changes according to their position with respect to the singular point and to the coordinate hyperplanes. These cells are shown to be contained in Bowen $(R,e^{-R})$-balls and we obtain an upper bound of $h(K)$ using an estimate on the number of such cells needed to cover $K$. We start with a description of the leaves of the foliation.  

For simplicity, we multiply $F$ with a constant in order to assume that
$$ \min \left\lbrace |\lambda_1|,\ldots,|\lambda_k| \right\rbrace=1 \quad \mbox{and} \quad 
\max \left\lbrace |\lambda_1|,\ldots,|\lambda_k| \right\rbrace=\lambda_*.$$
This does not change the foliation.
Write $\lambda_j=s_j+it_j$ with $s_j,t_j\in\R$.  
For $x=(x_1,\ldots,x_k)\in \D^k\setminus\{0\}$, define the holomorphic map $\varphi_x:\C\rightarrow\C^k\setminus\{0\}$ by
$$\varphi_x(\zeta):=\Big(x_1e^{\lambda_1\zeta},\ldots,x_ke^{\lambda_k\zeta}\Big)\quad \mbox{for}\quad \zeta\in\C.
$$
It is easy to see that $\varphi_x(\C)$ is the integral curve of $F$ which contains
$\varphi_x(0)=x$.

Write $\zeta=u+iv$ with $u,v\in\R$. The domain $\Pi_x:=\varphi_x^{-1}(\D^k)$ in $\C$ is defined by the inequalities
$$s_ju-t_jv< -\log|x_j| \quad \mbox{for}\quad j=1,\ldots,k.$$
So, $\Pi_x$ is a convex polygon which is not necessarily
bounded. It contains $0$ since $\varphi_x(0)=x$.  Moreover,  we have
$$
\dist(0,\partial \Pi_x)=\min \Big\lbrace -{\log|x_1|\over |\lambda_1|},\ldots,  -{\log|x_k|\over |\lambda_k|}  \Big\rbrace. 
$$
Thus, we obtain the following useful estimates
$$-\lambda_*^{-1}\log\|x\|_1\leq \dist(0,\partial \Pi_x)\leq -\log\|x\|_1.$$

The leaf of $\Lc$ through $x$ is given by $\widehat L_x:=\varphi_x(\Pi_x)$.
Observe that when the ratio $\lambda_i/\lambda_j$ are not all rational and all the coordinates of $x$ do not vanish, 
$\varphi_x:\Pi_x\to \widehat L_x$ is bijective and hence $\widehat L_x$ is simply connected. Otherwise, when the ratios $\lambda_i/\lambda_j$ are rational, all the leaves are closed submanifolds of $\D^k\setminus\{0\}$ and are biholomorphic to annuli.

Let $\tau_x:\D\to \Pi_x$ be a biholomorphic map such that $\tau_x(0)=0$. Then, $\widehat\phi_x:=\varphi_x\circ\tau_x$ is a map from $\D$ to $\widehat L_x$ and is a universal covering map of $\widehat L_x$ such that $\widehat\phi_x(0)=x$. If $\omega_P$ and $\omega_0$ denote the Hermitian forms associated to the Poincar\'e metric and the Euclidean metric on $\Pi_x$, define the function $\vartheta_x$ by
$$\omega_0=\vartheta_x^2\omega_P.$$
The following lemma describes the Poincar\'e metric on $\Pi_x$. The first assertion is probably known and is still valid if we replace $\Pi_x$ by an arbitrary convex domain in $\C$. Recall that $\widehat\eta$ is given by $\omega=\widehat\eta^2\omega_P$ on $\widehat L_x$. 

\begin{lemma}\label{lemma_vartheta}
We have for any $a\in\Pi_x$
$$ {1\over 2}\dist(a,\partial \Pi_x) \leq\vartheta_x(a)\leq  \dist(a,\partial \Pi_x).$$
In particular, we have
$$-{1\over 2}\lambda_*^{-1}\|x\|_1\log\|x\|_1\leq \widehat\eta(x)\leq -k\lambda_*\|x\|_1\log\|x\|_1.$$
\end{lemma}
\proof 
For each $a\in\Pi_x$, consider the family of holomorphic maps $\sigma:\D\to\Pi_x$ such that $\sigma(0)=a$. Denote by $\|D\sigma(0)\|$ the norm of the differential of $\sigma$ at 0 with respect to the Euclidean metrics on $\D$ and on  $\Pi_x$. The extremal property of $\omega_P$ implies that
$$\vartheta_x(a)={1\over 2}\max_\sigma\|D\sigma(0)\|.$$

Now, since the disc  with center $a$ and  radius $ \dist(a,\partial \Pi_x)$ is  contained in $\Pi_x$,
the  first estimate in the lemma follows.
Let $l$ be  a side  of $\Pi_x$ which is  tangent to  the above  disc. Consider the half-plane $\H$ containing $\Pi_x$ such that $l$ is contained in the boundary of $\H$. The Poincar\'e metric on $\Pi_x$ is larger than the one on $\H$. The Poincar\'e metric on $\H$ is associated to the Hermitian form $\dist(\cdot,\partial \H)^{-2}\omega_0$. This
implies the second estimate. 

For the second assertion in the lemma, we have
$$\widehat\eta(x)=\vartheta_x(0)\|D\varphi_x(0)\|.$$
This, the first assertion  and the above estimates on $\dist(a,\partial\Pi_x)$ 
imply the result because we get from the definition of $\varphi_x$ that 
$\|x\|_1\leq \|D\varphi_x(0)\|\leq k\lambda_*\|x\|_1$. 
\endproof

Fix a constant $\lambda$ such that $\lambda\geq\lambda_*$. For the main results in this section, it is enough to take $\lambda=\lambda_*$ but we will use later the case with a large constant $\lambda$.
Fix also a constant $0<\rho<1$ such that $K$ is strictly contained in $\rho\D^k$. 
Denote by $\Omega_x\subset \Pi_x$ the set of points $\zeta:=u+iv$ such that
$$\begin{cases}
s_ju-t_jv< -\log|x_j|-e^{-20\lambda R} \quad \mbox{for}\quad j=1,\ldots,k\\
|\zeta|\leq e^{20\lambda R}
\end{cases}$$
Recall from the introduction that 
$\alpha_1:=e^{-e^{7\lambda R}}$ and $\alpha_2:=e^{-e^{23\lambda R}}$.

\begin{lemma} \label{lemma_3R}
Assume that $\alpha_1\leq\|x\|_1\leq \rho$. Then, $\tau_x(\D_{7R})$ is contained in $\Omega_x$. 
In particular, $\widehat\phi_x(\D_{7R})$ is contained in $\rho'\D^k$ with $\rho':=e^{-e^{-21\lambda R}}\simeq 1-e^{-21\lambda R}$. 
\end{lemma}
\proof
It is not difficult to see that $\varphi_x(\Omega_x)$ is contained in $\rho'\D^k$. So,
the second assertion in the lemma is a direct consequence of the first one. We prove now the first assertion. 
Consider a point $\zeta\in \tau_x(\D_{7R})$. By the second inequality in Lemma \ref{lemma_vartheta}, the Poincar\'e distance between 0 and $\zeta$ in $\Pi_x$ is at least equal to
$$\int_0^{|\zeta|}{dt\over t+\dist(0,\partial\Pi_x)}\geq \int_0^{|\zeta|}{dt\over t-\log\|x\|_1}=\log\Big(1- {|\zeta|\over \log\|x\|_1}\Big)\cdot$$
Since $\zeta$ is a point in $\tau_x(\D_{7R})$, this distance is at most equal to $7R$. Therefore, using that $\|x\|_1\geq \alpha_1$, we obtain
$$|\zeta|\leq -\log\|x\|_1 e^{7R}\leq e^{20\lambda R}.$$

So, if the lemma were false, there would be a $\zeta$ with $\dist_P(0,\zeta)\leq 7R$ such that
$$s_ju-t_jv=-\log|x_j|-e^{-20\lambda R}$$
for some $j$.
It follows that $\dist(\zeta,\partial\Pi_x)\lesssim e^{-20\lambda R}$. Hence, using Lemma \ref{lemma_vartheta} and the estimate
$$\dist(0,\partial\Pi_x)\geq -\lambda^{-1}\log \|x\|_1\geq -\lambda^{-1}\log\rho,$$ 
we obtain
$$\dist_P(0,\zeta)\geq \Big|\int_{\dist(\zeta,\partial\Pi_x)}^{\dist(0,\partial\Pi_x)}{dt\over t}\Big|=\big|\log\dist(0,\partial\Pi_x) -\log \dist(\zeta,\partial\Pi_x)\big|>7R.$$
This is a contradiction.
\endproof

In order to prove Theorem \ref{th_local_entropy}, we need to study carefully the points which are close to a coordinate plane in $\C^k$. We have the following lemma. 

\begin{lemma} \label{lemma_epsilon_3}
Let $x$ be a point in $\rho\D^k$ such that $\|x\|_1\geq \alpha_1$. Let $m$ be an integer such that  $1\leq m\leq k$ and
$|x_j|\leq 2\alpha_2$ for $j=m+1,\ldots,k$. Then, $\widehat\phi_x(\D_{7R})$ is contained in $\D^m\times e^{-3R}\D^{k-m}$. 
\end{lemma}
\proof
Let $\xi$ be a point in $\D_{7R}$. Define $\zeta:=\tau_x(\xi)$ and $x':=\widehat\phi_x(\xi)=\varphi_x(\zeta)$. 
We have to show that $|x_j'|< e^{-3R}$ for $j\geq m+1$.
By definition of $\varphi_x$, we have 
$$|x'_j|=|x_j e^{\lambda_j \zeta}|\leq 2\alpha_2 e^{\lambda |\zeta|}.$$
This, combined with the first assertion of Lemma \ref{lemma_3R}, implies the result.
\endproof

We consider now the situation near the singular point.

\begin{lemma} \label{lemma_near_0}
If $\|x\|_1\leq 2\alpha_1$ and $\|y\|_1\leq 2\alpha_1$, then $x$ and $y$ are $(R,e^{-R})$-close.
\end{lemma}
\proof
It is enough to show that $\widehat\phi_x(\D_R)\subset e^{-2R}\D^k$. This and the similar property for $y$ imply the lemma. 
Consider $\xi$, $\zeta=u+iv$ and $x'$ as above with $\xi\in\D_R$. A computation as in the end of Lemma \ref{lemma_3R} implies that 
$$R\geq\dist_P(0,\zeta)\geq \log\dist(0,\partial\Pi_x) -\log \dist(\zeta,\partial\Pi_x).$$
Since $\|x\|_1\leq 2\alpha_1$, we have $\log\dist(0,\partial\Pi_x)> 6\lambda R$ and then $\dist(\zeta,\partial\Pi_x)\geq 3R$. It follows that 
$$|x_j'|=e^{\log |x_j|+s_ju-t_jv}\leq e^{-|\lambda_j|\dist(\zeta,\partial\Pi_x)}\leq e^{-2R}$$
for every $j$. The result follows.
\endproof

\begin{lemma} \label{lemma_epsilon_2}
Let $x$ be a point in $\rho\D^k$ and $1\leq m\leq k$ be an integer such that $\|x\|_1> \alpha_1$ and
$|x_j|\leq 2\alpha_2$ for $j=m+1,\ldots,k$. If $x':=(x_1,\ldots,x_m,0,\ldots,0)$, then $x$ and $x'$ are $(R,e^{-2R})$-close.
\end{lemma}
\proof
Fix a point $\xi\in\D_R$. We have to show that $\dist(\widehat\phi_x(\xi),\widehat\phi_{x'}(\xi))\leq e^{-2R}$.
Observe that by hypotheses $\Omega_{x'}=\Omega_x\subset \Pi_x$. 
So, by Lemma \ref{lemma_3R},
$\tau_{x'}$ defines a holomorphic map from $\D_{7R}$ to $\Pi_x$.
Observe that 
the Euclidean radius of $\D_{7R}$ is larger than $1-2e^{-7R}$. Hence, using the extremal property of the Poincar\'e metric,  we deduce that
$$\vartheta_{x'}(0)\leq (1+3e^{-7R}) \vartheta_x(0).$$  

Consider the map $\tau:=\tau_{x'}^{-1}\circ\tau_x$ from $\D_{7R}$ to $\D$. We have $\|D\tau(0)\|\geq 1-3e^{-7R}$. 
Composing $\tau_x$ with a suitable rotation allows us to assume that $D\tau(0)$ is a positive real number. By Lemma 2.3 in \cite{DinhNguyenSibony2} applied to $\tau$, there is a point $\xi'$ such that $\dist_P(\xi,\xi')\ll e^{-2R}$ and $\tau_{x'}(\xi')=\tau_x(\xi)$. 
Observe that the first $m$ coordinates of $\widehat\phi_x(\xi)$ are equal to the ones of $\widehat\phi_{x'}(\xi')$. Therefore,  by Lemma \ref{lemma_epsilon_3} applied to $x$ and to $x'$,  the distance between $\widehat\phi_x(\xi)$ and $\widehat\phi_{x'}(\xi')$ is less than $ke^{-3R}$. On the other hand,
$$\dist(\widehat\phi_{x'}(\xi),\widehat\phi_{x'}(\xi'))\lesssim \dist_P(\widehat\phi_{x'}(\xi),\widehat\phi_{x'}(\xi'))=\dist_P(\xi,\xi')\ll e^{-2R}.$$
The lemma follows.
\endproof

Consider now two points $x$ and $y$ in $\rho\D^k$ such that for each $1\leq j\leq k$ one of the following properties holds
\begin{itemize}
\item[(S1)]  $|x_j|<\alpha_2$ and $|y_j|<\alpha_2;$  
\item[(S2)] $x_j,y_j\not=0$ and   $\big|{ x_j\over y_j} -1\big |< e^{-22\lambda R}$ and $\big|{ y_j\over x_j} -1\big |<e^{-22\lambda R}.$
\end{itemize}

We have the following proposition.

\begin{proposition} \label{prop_close_W}
Under the above conditions, $x$ and $y$ are $(R,e^{-R})$-close.
\end{proposition}
\proof
If $\|x\|_1\leq 2\alpha_1$ and $\|y\|_1\leq 2\alpha_1$, then Lemma \ref{lemma_near_0} implies the result. Assume this is not the case. By condition (S2), we have $\|x\|_1\geq\alpha_1$ and $\|y\|_1\geq \alpha_1$. Moreover, up to a permutation of coordinates, we can find $1\leq m\leq k$ such that $|x_j|\geq \alpha_2$, $|y_j|\geq \alpha_2$ for $j\leq m$ and $|x_j|\leq 2\alpha_2$, $|y_j|\leq 2\alpha_2$ for $j\geq m+1$. By Lemma \ref{lemma_epsilon_2}, we can assume that $x_j=y_j=0$ for $j\geq m+1$. Now, in order to simplify the notation, we can assume without loss of generality that $m=k$. So, we have $|x_j|\geq\alpha_2$ and $|y_j|\geq \alpha_2$ for every $j$. 

Define the linear holomorphic map $\Psi_{x,y}:\C^k\to\C^k$ by 
$$\Psi_{x,y}(z):=\Big({y_1\over x_1} z_1,\ldots,{y_k\over x_k} z_k\Big),\qquad  z=(z_1,\ldots,z_k)\in\C^k.$$
This map preserves the foliation and sends $x$ to $y$. Moreover, the property (S2) implies that for such $x,y$,  
$\|\Psi_{x,y}-\id\|\ll e^{-21\lambda R}$ on $\D^k$. 

Define also $\widetilde\phi_y:=\Psi_{x,y}\circ\widehat\phi_x$. It follows from the last assertion in Lemma \ref{lemma_3R} that this map is well-defined on $\D_{7R}$ with image in  $\widehat L_y$ and we have $\widetilde\phi_y(0)=y$. 
We also deduce from (S2) that $\dist(\widehat\phi_x(\xi),\widetilde\phi_y(\xi))\ll e^{-R}$. 
It remains to compare $\widehat\phi_y$ and $\widetilde\phi_y$ on $\D_R$.

Using the extremal property of the Poincar\'e metric, we obtain
$$\|D\widehat\phi_y(0)\|\geq (1-e^{-6R})\|D\widetilde\phi_y(0)\|\geq (1-e^{-5R})\|D\widehat\phi_x(0)\|.$$
By symmetry, we deduce that $\|D\widehat\phi_x(0)\|$, $\|D\widehat\phi_y(0)\|$ and $\|D\widetilde\phi_y(0)\|$ are close, i.e. their ratios are bounded by $1+e^{-4R}$. 

Since $\widehat\phi_y$ is a universal covering map, there is a unique holomorphic map $\tau:\D_{7R}\to\D$ such that $\tau(0)=0$ and  $\widehat\phi_y\circ \tau= \widetilde\phi_y$. Composing $\widehat\phi_y$ with a suitable rotation allows us to assume that $D\tau(0)$ is a positive real number. It follows from the above discussion that $|1-\tau'(0)|\leq e^{-4R}$. Lemma 2.3 in \cite{DinhNguyenSibony2} implies that $\dist_P(\xi,\tau(\xi))\ll e^{-R}$ on $\D_R$. Hence, $\dist_P(\widehat\phi_y(\xi),\widetilde\phi_y(\xi))\ll e^{-R}$ for $\xi$ in $\D_R$. The proposition follows. 
\endproof 

\noindent
{\bf End of the proof of Theorem \ref{th_local_entropy}.}
The idea is to divide $\D^k$ into cells which are contained in Bowen $(R,e^{-R})$-balls. 
The sizes of these cells are very different and this is one of the main difficulties in the proof of Theorem \ref{th_main}. 

We first divide $\D$ into rings using the circle of center 0 and of radius $\alpha_2e^{ne^{-23\lambda R}}$ for $n=1,\ldots, e^{46\lambda R}$. In fact, we have to take the integer part of the last number but we will not write it in order to simplify the notation.
Then, we divide these rings into cells using $e^{23\lambda R}$ half-lines starting at  0 which are equidistributed in $\C$. 
We obtain less than $e^{70\lambda R}$ cells and we denote by $\Sigma$ the set of the vertices, i.e. the intersection of circles and half-lines. 
Except those at 0, if a cell contains a point $a$, it looks like a rectangle whose sides are approximatively $|a|e^{-23\lambda R}$.
Consider the product of $k$ copies of $\D$ together with the above division, 
we obtain a division of $\D^k$ into less than $e^{70\lambda k R}$ cells.

Consider two points $x,y$ in $\rho\D^k$ which belong
 to the same cell. They satisfy the conditions (S1) and (S2). So, by Proposition \ref{prop_close_W}, they are $(R,e^{-R})$-close. It follows that if a cell is contained in $\rho\D^k$, it is contained in a Bowen  $(R,e^{-R})$-ball. We deduce that the entropy of $K$ is bounded by $70\lambda k$. The estimate holds for $\lambda=\lambda_*$.
\hfill $\square$

\begin{proposition} The function $\widehat\eta$ is locally H\"older continuous outside the coordinate hyperplanes $\{x_j=0\}$, $1\leq j\leq k$, with H\"older exponent $(6\lambda_*)^{-1}$. 
\end{proposition}
\proof
We consider $x$ in a fixed compact outside the coordinate hyperplanes. So, as in the case without singularities \cite{DinhNguyenSibony2}, we can show that $\widehat\eta(x)$ is bounded from above and from below by strictly positive constants.
Using the comparison between $D\widehat\phi_x(0)$ and $D\widehat\phi_y(0)$ in the proof of Proposition \ref{prop_close_W}, we deduce that 
if $\dist(x,y)\leq e^{-23\lambda R}$, then $|\widehat\eta(x)-\widehat\eta(y)|\leq e^{-4 R}$. The result follows.
\endproof

Note that the above division of $\D^k$ respects the invariance of the foliation under the homotheties. However, it is important to observe that the Poincar\'e metric on leaves is not invariant under the homotheties, see Lemma \ref{lemma_vartheta}. As a consequence, 
the plaques in the above cells are very small in the sense that their Poincar\'e diameters tend to 0 when $R$ tends to infinity. We will heavily use properties of $(\D^k,\Lc,\{0\})$, in particular the above division into cells, as a model of singular flow boxes in our study of global foliations. We give now some construction that will be used later.

Recall that $\Sigma$ is the intersection of the circles and the half-lines used in the end of the proof of Theorem \ref{th_local_entropy}. Note that this set depends on $R$. 
We introduce now $p$ maps $J_l$ from $\Sigma^k\cap {3\over 4}\D^k$ to $\Sigma^k$, 
where $p$ is a large integer that will be fixed just before Lemma \ref{lemma_p_points_1} below.
These maps describe roughly, the displacement of points when we travel along leaves following $p$ given directions after a certain time, which is independent of $R$. 

For $x\in\Sigma^k$ and $0\leq l\leq p-1$, define 
$$z^l:=\varphi_x\Big(-t\log\|x\|_1e^{2i\pi l\over p}\Big),$$
where the constant $0<t\ll \lambda_*^{-1}$ will be fixed just after Lemma \ref{lemma_p_points_1} below.  
For the moment, we will need that the real part of $te^{2i\pi l\over p}\lambda_j$ is not equal to $-1$ for all $j$ and $l$. This property simplifies the proof of Lemma \ref{lemma_corr} below.
We choose $J_l(x)$ a vertex of a cell of $\D^k$ which contains $z^l$. Note that the choice is not unique but this is not important for our problem.

The point $x$ is displaced to the points $z^l$ when we travel 
following $p$ directions after a certain time. The points $J_l(x)$ give us an approximation of $z^l$ and allow us to understand the displacement of points near $x$ using the lattice $\Sigma^k$. We will use these maps $J_l$ in order to make an induction on hyperbolic time $R$ which is an important step in the proof of Theorem \ref{th_main}, see also Lemmas \ref{lemma_cover_D} and \ref{lemma_cover_D_bis} below.

\begin{lemma} \label{lemma_corr}
Let $y=(y_1,\ldots,y_k)$ be a point in $\Sigma^k$ and $0\leq l\leq p-1$ be an integer. There is a constant $M>0$ independent of $R$ such that if $|y_j|>\alpha_1$ for every $j$, then $J_l^{-1}(y)$ contains at most $M$ points. 
\end{lemma}
\proof
Let $x=(x_1,\ldots,x_k)$ be such that $J_l(x)=y$. For simplicity, write  $z$ for the point $z^l=\varphi_x\big(-t\log\|x\|_1e^{2i\pi l\over p}\big)$. Since $z$ and $y$ belong to the same cell, $z_j/y_j$ is very close to 1. We deduce that $|z_j|\geq {1\over 2}\alpha_1$ for every $j$. 
Using the definition of $\varphi_x$, we obtain 
$${1\over 2}\alpha_1\leq |z_j|=|x_j|\big| e^{-t\log\|x\|_1 e^{2i\pi l/p}}\big|\leq |x_j| \|x\|_1^{-t}\leq |x_j|^{1-t}.$$
Hence, since $t$ is small, we infer $|x_j|\geq \alpha_2$ for every $j$.

We will only consider the case where $\|x\|_1=|x_1|$. The other cases are treated in the same way. In order to simplify the notation, assume also that $l=0$. 
Consider another point $x'$ satisfying similar properties. It is enough to show that the number of such points $x'$ is bounded. Define $z':=\varphi_{x'}\big(-t\log|x'_1|\big)$.

We have
$$z=(x_1|x_1|^{s\lambda_1},x_2|x_1|^{s\lambda_2},\ldots, x_k|x_1|^{s\lambda_k})$$
and 
$$z'=(x'_1|x'_1|^{s\lambda_1},x'_2|x'_1|^{s\lambda_2},\ldots, x'_k|x'_1|^{s\lambda_k}).$$
Since $z,y$ belong to the same cell and $z',y$ satisfy the same property, we have $|z_j'/z_j-1|\lesssim e^{-23\lambda R}$ and $|z_j/z_j'-1|\lesssim e^{-23\lambda R}$. 
Recall that $|x_j|\geq \alpha_2$ and $|x_j'|\geq \alpha_2$ and 
$x_j,x_j'$ are in $\Sigma$. Using that the real part of $t\lambda_j$ is not equal to $-1$, we easily see that for each $j$ there is a bounded number of $x'_j$ satisfying the above properties. For this purpose, we can also use a homothety in order to reduce the problem to the case where $|x_j|\simeq 1$. 
\endproof

In general, the point $J_l(x)$ does not belong to a plaque containing  $z^l$ but it is however very close to such a plaque. Write $J_l(x)=w=(w_1,\ldots,w_k)$. We have the following lemma.

\begin{lemma} \label{lemma_J_l_close}
Assume that $|z^l_j| >\alpha_2$ for every $j$. Then, the plaque $L_{z^l}[\|z^l\|_1\epsilon_0]$ intersects $\{w_1\}\times \D^{k-1}$ at a unique point $w^l$. Moreover, there is a constant $\gamma\geq 1$ independent of $\lambda$ and $R$ such that $|w_j/w^l_j-1|\leq \gamma e^{-23\lambda R}$ and 
$|z^l_j/w^l_j-1|\leq \gamma e^{-23\lambda R}$ for every $j$.
\end{lemma}
\proof
We show that the lemma is true for a general point $z$ with $|z_j| >\alpha_2$ for every $j$ and for any vertex $w$ of a cell containing $z$. Observe that $w_j\not=0$ for every $j$. 
We first consider the case where $|z_j|\simeq 1/2$ for every $j$. In this case, $L_z[\epsilon_0/2]$ is the graph of a map with bounded derivatives over a domain on the first coordinate axis. The cell containing $z$ looks like a cube of size $\simeq e^{-23\lambda R}$. So, the lemma is clear in this case. 

Consider now the general case. Observe that the foliation and the set $\Sigma^k$ are invariant when we multiply a coordinate by a power of $e^{e^{-23\lambda R}}$. Therefore, multiplying the coordinates by a same constant allows us to assume that $\|z\|_1=|z_j|\simeq 1/2$ for some index $j$. Then, $L_z[\epsilon_0/2]$ is a graph of a map with bounded derivatives over a domain $D$ in the $j$-th axis.
We claim that if we multiply each coordinate $z_l$ with $l\not=j$ with an appropriate power of $e^{e^{-23\lambda R}}$, the problem is reduced to the first case.
Indeed, the image of $L_z[\epsilon_0/2]$  is still  the graph over $D$ of a map with bounded derivatives because this graph is contained in a plaque. Therefore, the size of $L_z[\epsilon_0/2]$ changes with a factor bounded independently of $R$ and $\lambda$. Since $\epsilon_0$ is small enough, the same arguments as in the first case give the result.
\endproof

Observe that in the proof of Theorem \ref{th_local_entropy}, the holomorphic maps $\Psi_{x,y}$ are used  in order to control the distance between the leaves. Recall that for global foliations without singularities, we have to use some orthogonal projections in order to construct a parametrization $\psi$ as in the definition of conformally $(R,\delta)$-close points. These maps 
are not holomorphic but we can correct them using Beltrami's equation. In order to prove Theorem \ref{th_main}  for the case with singularities, we will use both kinds of maps depending if we are near or far from singular points. We describe below some relations between these two kinds of maps and show that one can glue them near the boundaries of singular flow boxes.

Recall that for a constant $\epsilon_0$ small enough, if $x,y$ are two points in ${3\over 4}\D^k\setminus {1\over 4}\D^k$ such that $\dist(x,y)\leq \epsilon_0$, then the basic projection 
$\Phi$ from $L_x[\epsilon_0]$ to $L_y[10\epsilon_0]$  is well-defined with image in $L_y[3\epsilon_0]$. Fix a constant 
$\kappa>1$ large enough such that the following estimate holds
on $L_x[\epsilon_0]$:
$$ \|\Phi-\id\|_{\Cc^{2}}\leq  \kappa  \dist (\Phi(x),x),$$
where we compute the norm using the Euclidean metric on $\widehat L_x$ and $\widehat L_y$
 
The situation is more delicate near the singularities. We use the fact that the foliation in $\D^k$ is invariant under homotheties $z\mapsto t z$ with $|t|\leq 1$. Such a homothety multiplies the distance with $|t|$.  For $x\in {3\over 4}\D^k$ and $t:={4\over 3}\|x\|_1$, some point in the boundary of ${3\over 4}\D^k$ is sent to $x$. Thus,  if 
$\dist(x,y)\leq \|x\|_1\epsilon_0$, we can define the orthogonal projection $\Phi$ from  $L_x[\|x\|_1\epsilon_0]$ to  $L_y[10\|x\|_1\epsilon_0]$ with image in $L_y[3\|x\|_1\epsilon_0]$. Moreover, we have
$$ \|\Phi-\id\|_{\Cc^0}\leq  \kappa \dist (\Phi(x),x)
\quad \mbox{and}\quad \|\Phi-\id\|_{\Cc^2}\leq  \kappa \dist (\Phi(x),x) \|x\|_1^{-2}.$$
Note that $L_x[\|x\|_1\epsilon_0]$ and  $L_y[3\|x\|_1\epsilon_0]$ are connected and simply connected.

We see that  when we are near the singularity, the Beltrami coefficient of $\Phi$ may have a large $\Cc^1$-norm and cannot satisfy the estimates in the definition of conformally $(R,\delta)$-close points. 
This is the reason why we have to use the holomorphic maps $\Psi_{x,y}$ defined above. Nevertheless, the control of the $\Cc^0$-norm $\|\Phi-\id\|_{\Cc^0}$ is still good near the singular point. Therefore, in the proof of Theorem \ref{th_main}, it is convenient to use first the orthogonal projections and then correct the parametrization using the two lemmas below.

\begin{lemma} \label{lemma_proj_E}
Let $x,y$ be two points in ${1\over 2}\D^k$ such that $y_j\not=0$ and  $|x_j/y_j-1|\leq e^{-10R}$ for every $j$. Then, there is an orthogonal projection $\Phi_{x,y}$ from $\widehat L_x\cap {3\over 4}\D^k$ to $\widehat L_y$ which coincides on $L_z[\|z\|_1\epsilon_0]$ with the basic projection associated to $z$ and to $w:=\Psi_{x,y}(z)$ for every $z\in \widehat L_x\cap{3\over 4}\D^k$. Moreover, we have
$\|\Phi_{x,y}-\id\|_{\Cc^0}\leq e^{-9R}$. 
\end{lemma}
\proof
Observe that 
$$|z_j/w_j-1|=|x_j/y_j-1|\leq e^{-10R}.$$
Therefore, there is a basic projection from $L_z[\|z\|_1\epsilon_0]$ to $L_w[10\|z\|_1\epsilon_0]$ with image in $L_z[3\|z\|_1\epsilon_0]$ described above. It satisfies the estimate in the lemma. 
It is not difficult to see that such projections coincide on the intersection of their domains of definition. 
Indeed, the property is clear near the boundary of ${3\over 4}\D^k$. The general case can be reduced to that case using a homothety.
The lemma follows.
\endproof

\begin{lemma} \label{lemma_proj_hol}
Let $x,y$ be as in Lemma \ref{lemma_proj_E}. Then, there is a map $\widetilde\Psi_{x,y}$ from $\widehat L_x\cap {3\over 4}\D^k$ to $\widehat L_y$ such that $\widetilde\Psi_{x,y}=\Psi_{x,y}$ on $\widehat L_x\cap {1\over 4}\D^k$ and  
$\widetilde\Psi_{x,y}=\Phi_{x,y}$ on $\widehat L_x\cap {3\over 4}\D^k\setminus {1\over 2}\D^k$. 
In particular, $\widetilde\Psi_{x,y}$ is holomorphic in ${1\over 4}\D^k$. 
Moreover, we have $\|\widetilde\Psi_{x,y}-\id\|_{\Cc^0}\leq e^{-8R}$ on $\widehat L_x\cap {3\over 4}\D^k$ and $\|\widetilde\Psi_{x,y}-\id\|_{\Cc^2}\leq e^{-8R}$ on $\widehat L_x\cap {3\over 4}\D^k\setminus {1\over 4}\D^k$. 
\end{lemma}
\proof
Define $z':=\Phi_{x,y}(z)$ and  $w:=\Psi_{x,y}(z)$.
So, $z'$ belongs to a small plaque containing $w$. Using a homothety and the definition of $\varphi_w$, we can check that $|z_i'/w_i-1|\lesssim e^{-10R}$. Therefore, 
the function $\log(z'_i/w_i)$ is well-defined. 

Let $0\leq \chi \leq 1$ be a smooth function equal to 0 on ${1\over 4}\D^k$ and equal to 1 outside ${1\over 2}\D^k$. 
Define $\chi_i(z):=\chi(z) \log(z'_i/w_i)$ and 
$$\widetilde\Psi_{x,y}(z):=\Big(e^{\chi_1(z)}w_1,\ldots, e^{\chi_k(z)}w_k\Big).$$
Clearly, $\widetilde\Psi_{x,y}=\Psi_{x,y}$ on $\widehat L_x\cap {1\over 4}\D^k$ and  
$\widetilde\Psi_{x,y}=\Phi_{x,y}$ on $\widehat L_x\cap {3\over 4}\D^k\setminus {1\over 2}\D^k$. 

Observe that a point $w'=(w'_1,\ldots,w'_k)$ belongs to a small plaque containing $w$ if and only if 
$$\lambda_j\log(w'_i/w_i)=\lambda_i \log(w'_j/w_j).$$
This holds in particular for $w'=z'$ and hence, $\widetilde\Psi_{x,y}(z)$ satisfies also this criterium. It follows that $\widetilde\Psi_{x,y}$ has values in $\widehat L_y$. 

We deduce from the definition of $\Psi_{x,y}$ that  $\|\Psi_{x,y}-\id\|_{\Cc^0}\leq e^{-9R}$ on $\widehat L_x\cap {3\over 4}\D^k$ and $\|\Psi_{x,y}-\id\|_{\Cc^2}\leq e^{-9R}$ on $\widehat L_x\cap {3\over 4}\D^k\setminus {1\over 4}\D^k$. Recall that $\Phi_{x,y}$ satisfies similar estimates. So, $\widetilde\Psi_{x,y}$ satisfies the estimates in the lemma. 
\endproof


\section{Poincar\'e metric on leaves} \label{section_poincare}

In this section, we will meet another important difficulty for our study: a leaf of the foliation may visit singular flow boxes without any obvious rule. However, we analyze the behavior and get an explicit estimate on the modulus of continuity of the Poincar\'e metric on leaves.
We are concerned with the following class of foliations. 
 
\begin{definition}\label{uniform_hyperbolic_laminations} \rm
 A   Riemann surface foliation with singularities $(X,\Lc,E)$  on a Hermitian compact complex manifold $X$
 is  said  to be {\it Brody hyperbolic}  if  there is  a constant  $c_0>0$  such that
$$\|D\phi(0)\| \leq  c_0$$ 
for all  holomorphic maps  $\phi$ from $\D$  into  a leaf.
\end{definition}

It is clear that if the foliation is Brody hyperbolic then its leaves are hyperbolic in the sense of Kobayashi.
Conversely, the  Brody hyperbolicity is  a consequence of
the non-existence  of  holomorphic non-constant maps
$\C\rightarrow X$  such that out of  $E$ the image of $\C$ is  locally contained in leaves, 
see \cite[Theorem 15]{FornaessSibony2}. 
 
On the other hand,   Lins Neto proved  in   \cite{Neto} that  for every   holomorphic foliation of degree larger than 1  in $\P^k$, with non-degenerate singularities,  there is a smooth metric with negative curvature on its tangent bundle, see also Glutsyuk \cite{Glutsyuk}.
Hence, these foliations are Brody hyperbolic.
Consequently,  holomorphic foliations  in $\P^k$  are generically  Brody hyperbolic, see also  \cite{NetoSoares}.

From now on, we assume that $(X,\Lc,E)$ is Brody hyperbolic.
We have the following result with the notation introduced above.

\begin{theorem} \label{th_poincare_C0}
Let $(X,\Lc,E)$ be a Brody hyperbolic foliation by Riemann surfaces on a Hermitian compact complex manifold $X$. Assume that the singular set $E$  is  finite  and  that all  points of $E$  are  linearizable. 
 Then, there are constants  $c>0$ and $0<\alpha<1$   such that
 $$|\eta(x)-\eta(y)|\leq  c\Big ( {\max\{\lof \dist(x,E), \lof\dist(y,E)\}\over \lof\dist(x,y)}\Big)^\alpha$$
 for all $x,y$ in $X\setminus E$. 
\end{theorem}

Note that in \cite{Shcherbakov}, Shcherbakov proved a surprising result which says that when $X$ is a projective space,  $\eta$ is smooth on $\P^k\setminus E$. However, this result does not implies Theorem \ref{th_poincare_C0} for the case of projective spaces. We thank Marco Brunella for this reference.

The following estimates are crucial in our study. 

\begin{proposition} \label{prop_eta_E}
Under the hypotheses
of Theorem  \ref{th_poincare_C0},  there exists  a  constant $c_1>1$  such that $\eta \leq c_1$ on $X$, $\eta\geq c_1^{-1}$ outside the singular flow boxes ${1\over 4}\U_e$ and 
$$c_1^{-1} s \lof s \leq\eta(x)  \leq c_1  s \lof s$$
for $x\in X\setminus E$  and $s:=\dist(x,E)$.  
\end{proposition}
\proof 
It is enough to prove the last assertion. 
Without loss of generality, we only need to show it for
all $x$  in a singular flow box ${1\over 2}\U_e$. We identify $\U_e$ with the model in Section \ref{section_local} and use the notation introduced there. Recall that for simplicity, we use here a metric on $X$ whose restriction to $\U_e\simeq\D^k$ coincides with the Euclidean metric. 

Define the map $\tau:\D\to L_x$ by 
$$\tau(\xi):=\varphi_x(-c^{-1}\log\|x\|_1\xi)$$
for some constant $c>1$ large enough. Since $\|D\varphi_x(0)\|$ is bounded from below by $\|x\|_1$, we deduce that 
$$\eta(x)\geq \|D\tau(0)\|\gtrsim -\log \|x\|_1\|x\|_1 \gtrsim  s\lof s.$$
This gives us the first inequality in the last assertion of the proposition.

Next, since the foliation is Brody  hyperbolic, the function $\eta$ is bounded from above. It follows that there  is  a  constant $R_0>0$ independent of $x\in {1\over 2}\U_e$ such that  
the   disc of center $x$ in $L_x$ with radius $R_0$ with respect to
the Poincar\'e  metric is  contained  in $\U_e$.
Let $r_0\in (0,1)$ such that  $\D_{R_0}=r_0\D$.   
We have $\phi_x(r_0\D)\subset \D^k$. We then deduce from the extremal property of the Poincar\'e metric that $\eta(x)\leq r_0^{-1}\widehat\eta(x)$.
Lemma \ref{lemma_vartheta} implies the result.
 \endproof 

The following lemma gives us a speed estimate when we travel in a singular flow box along a geodesic with respect to the Poincar\'e metric on leaves. We denote by $[0,\xi]$ the  segment joining $0$ and $\xi$ in $\D$.
 
\begin{lemma}\label{lemma_close_sing}
There is a constant $c_2>0$ independent of $x$ in ${1\over 2}\U_e\simeq {1\over 2}\D^k$  with the following property. If $\xi\in \D_R$ such that $\phi_x([0,\xi])\subset {1\over 2}\D^k$ and if $y:=\phi_x(\xi)$, then
$$\|x\|_1^{e^{c_2R}} \leq\|y\|_1  \leq \|x\|_1^{e^{-c_2R}}.$$ 
\end{lemma}
\proof 
We only have to prove the first inequality. The second one is obtained by exchanging  $x$ and $y$. So, we only need to consider the case where $\|y\|_1\leq \|x\|_1$. 
By Proposition \ref{prop_eta_E}, we have
$$R\geq \dist_P(x,y)\gtrsim \int_{\|y\|_1}^{\|x\|_1} {dt\over t|\log t|} =\log {\log\|y\|_1\over \log\|x\|_1}\cdot$$
The lemma follows.
\endproof

The following lemma shows us how deep a leaf can go into a singular flow box before the hyperbolic time $R$.

\begin{lemma} \label{lemma_close_sing_bis}
There is a constant $c_3>0$ such that
for every $x\in X\setminus E$, we have
$$ \dist(\phi_x(\D_R),E)\geq e^{-\lof\dist(x,E)e^{c_3R}}.$$
\end{lemma}
\proof
We only have to estimate $\dist(\phi_x(\xi),E)$ for a point $\xi\in\D_R$ such that $\phi_x(\xi)$ is close to a singular point $e$. So, we can replace $x$ by a suitable point in $\phi_x([0,\xi])$ in order to assume that $\phi_x([0,\xi])\subset {1\over 2}\U_e$. Consequently, Lemma \ref{lemma_close_sing}  implies the result.
\endproof
  
Recall that we define the notion of conformally $(R,\delta)$-close as in 
\cite{DinhNguyenSibony2} using the constant $A:=3c_1^2$ and a number $0<\delta\leq e^{-2R}$. In order to prove Theorem \ref{th_poincare_C0}, we follow the approach of that paper. 
We will need the following result for a large enough constant $c_4$.

\begin{proposition} \label{prop_quasi_close}
Let $c_4>1$ be a fixed constant.
Let $x$ and $y$ be conformally $(R,\delta)$-close such that $c_4^{-1}\eta(y)\leq \eta(x)\leq c_4\eta(y)$ and $e^{2R}\delta\leq \eta(y)$.  
Then, there is a real number $\theta$ such that if $\phi_y'(\xi):=\phi_y(e^{i\theta}\xi)$, we have
$$|\eta(x)-\eta(y)|\leq e^{-R/4} \quad \mbox{and}\quad 
\dist_{\D_{R/3}} (\phi_x,\phi_y')\leq e^{-R/4}.$$
\end{proposition}
\begin{proof}
We argue  as in the  proof of Proposition 2.2 in \cite{DinhNguyenSibony2} with the same notation. The difference with the non-singular case is that we have a
condition on the relative size of $\eta(x)$ and $\eta(y)$.
We get that
$$\|D\phi_y(0)\| - \|D\widetilde\phi_y(0)\| \lesssim e^R\delta +e^{-R}(\eta(x)+\eta(y)).$$
Recall that $\eta(y)=\|D\phi_y(0)\|$. We deduce from the hypotheses that the constant $\lambda$ used in the proof of 
 Proposition
2.2  in \cite{DinhNguyenSibony2} satisfies
$$1-\lambda \lesssim  {e^R\delta +e^{-R}(\eta(x)+\eta(y))  \over \eta(y) }\lesssim e^{-R}.$$ 
Now, it is enough to follow the proof of the above proposition. We use here that $e^{-R/3}\ll e^{-R/4}$ since we only consider $R$ large enough.
\end{proof}
 
In order to apply the last proposition, we have to show that if $x$ and $y$ are close enough, they are conformally $(R,\delta)$-close. So, we need to construct a map $\psi$ as in the definition of conformally $(R,\delta)$-close points. As in the case without singularities, $\psi$ will be obtained by composing $\phi_x$ with basic projections from $L_x$ to $L_y$. 
There are two main steps. The first one is to show that up to time $R$, the leaves $L_x$ and $L_y$ are still close enough in order to define basic projections. The second one is that we can glue these projections together in order to get a well-defined map on $\D_R$. The second step can be treated as in the case without singularities. So, for simplicity,  in what follows, we only consider the first problem.

Recall that if $x,y$ are two points outside the singular flow boxes ${1\over 2}\U_e$ such that $\dist(x,y)\leq \epsilon_0$, then the basic projection 
$\Phi$ from  $L_x[\epsilon_0]$ to $L_y[3\epsilon_0]$ is well-defined. Moreover, we have the following estimate on $L_x[\epsilon_0]$ for a fixed constant $\kappa>1$:
$$ \|\Phi-\id\|_{\Cc^{2}}\leq  \kappa  \dist (\Phi(x),x).$$

Inside the singular flow boxes, the foliation and the metric are invariant under homotheties. Therefore, since $\epsilon_0$ is small enough, we deduce that for all $x,y\in X$ such that
$\dist(x,y)\leq \dist(x,E)\epsilon_0$, there is a basic projection $\Phi$ from $L_x[\dist(x,E)\epsilon_0]$ to $L_y[3\dist(x,E)\epsilon_0]$ which satisfies
$$ \|\Phi-\id\|_{\Cc^{2}}\leq  \kappa {\dist (\Phi(x),x)\over \dist(x,E)^2}\cdot$$
It is not difficult to obtain the following useful estimates where we change the constant $\kappa$ if necessary:
$${\dist (\Phi(z),z)\over \dist(z,E)^2}      \leq  \kappa  {\dist (\Phi(x),x)\over \dist(x,E)^2}\quad \mbox{and}\quad {\dist (\Phi(z),z)\over \dist(z,E)^6}      \leq  \kappa  {\dist (\Phi(x),x)\over \dist(x,E)^6}$$
for every $z$ in $L_x[\dist(x,E)\epsilon_0]$.

\begin{proposition}\label{prop_conformal_close}
There exists a constant $c_5>1$  with the  following property. Let  $x,y\in X\setminus E$, $R$ large enough and 
$\delta:=  \dist(x,y)e^{\lof \dist(x,E)e^{c_5R}}$ such that $\delta\leq e^{-2R}$.
Then,  $x$ and $y$ are conformally $(R,\delta)$-close.  
\end{proposition}

Taking for granted this proposition, we complete  the proof of Theorem \ref{th_poincare_C0}.

\medskip
\noindent
{\bf End of the proof of Theorem \ref{th_poincare_C0}.}
We can assume that $x,y$ are close and 
$$|\log\dist(x,y)|\gg \lof\dist(x,E) \quad \mbox{and} \quad |\log\dist(x,y)|\gg \lof\dist(y,E).$$
In particular, we have 
$$\dist(x,y)\ll \dist(x,E) \quad \mbox{and}\quad \dist(x,y)\ll \dist(y,E).$$ 
We will apply Propositions \ref{prop_quasi_close} and \ref{prop_conformal_close}. 
Choose $R>1$  such that
$$|\log \dist(x,y)|=\lof\dist(x,E)e^{2c_5R}.$$
So, $R$ is a large number. 
By Proposition \ref{prop_conformal_close}, $x$ and $y$ are conformally $(R,\delta)$-close with
$$\delta:= \dist(x,y) e^{\lof\dist(x,E)e^{c_5R}}=e^{\lof\dist(x,E)(e^{c_5R}-e^{2c_5 R})}.$$

It is clear that $\delta\leq e^{-2R}$. The above identities, together with 
Proposition \ref{prop_eta_E}, also imply that $e^{2R}\delta\leq \eta(y)$. 
Therefore,  $x,y,R$ and $\delta$ satisfy the hypotheses of Proposition \ref{prop_quasi_close}. Consequently,   
$$|\eta(x)-\eta(y)|\leq e^{-R/4}=\Big({\lof\dist(x,E)\over |\log\dist(x,y)|}\Big)^{1/(8c_5)}.$$
The result follows.
\hfill $\square$ 
   
\bigskip
   
We prove now Proposition \ref{prop_conformal_close}. Consider $x,y,R$ and $\delta$ as in this proposition. So,  $x,y$ are close and  satisfy 
$$|\log \dist(x,y)|\geq \lof\dist(x,E)e^{c_5R}\quad \mbox{and}\quad |\log\dist(x,y)|\geq \lof\dist(y,E)e^{c_5R}.$$
We have to construct a map $\psi:\overline\D_R\to L_y$ as in the definition of conformally $(R,\delta)$-close points.
We first consider the case where $y$ is the orthogonal projection of $x$ to $L_y[3\dist(x,E)\epsilon_0]$.

Fix a point $\xi\in\overline\D_R$. We will construct the map $\psi$ on a neighbourhood of $[0,\xi]$ using basic projections from $L_x$ to $L_y$. This allows us to prove that $\psi$ is well-defined on $\overline\D_R$ as in the case without singularities.  The details for  this last point are the same as in the non-singular case.

We first divide $[0,\xi]$ into a finite number of segments $[\xi^j,\xi^{j+1}]$ with $0\leq j\leq n-1$ and 
define $x^j:=\phi_x(\xi^j)$. The points $\xi^j$ are chosen by induction: 
$\xi^0:=0$ and  $\xi^{j+1}$ is the closest point to $\xi^j$ satisfying
$$\dist(x^{j+1},x^j)=\dist(x^j,E)\epsilon_1$$
with a fixed constant $\epsilon_1\ll\epsilon_0$.
The integer $n$ satisfies $\xi\in [\xi^{n-1},\xi^n]$. 

Define also by induction the points $y^j$ in $L_y$ with $y^0:=y$ and $y^{j+1}$ is the image of $x^{j+1}$ by the basic  projection $\Phi_j:L_{x^j}[\dist(x^j,E)\epsilon_0]\to L_{y^j}[3\dist(x^j,E)\epsilon_0]$. Note that since $\epsilon_1\ll\epsilon_0$, the point $y^{j+1}$ is also the image of $x^{j+1}$ by $\Phi_{j+1}$. We deduce from the properties of basic projections that
$${\dist(y^{j+1},x^{j+1})\over \dist(x^{j+1}, E)^6} \leq \kappa {\dist(y^j,x^j)\over \dist(x^j, E)^6}\cdot$$    
The following lemma guarantees, by induction on $j$, the existence of the above projections $\Phi_j$.

\begin{lemma} \label{lemma_delta}
We have 
$$\dist(y^j,x^j)\leq e^{-2R} \delta \dist(x^j,E)^6$$
for $j=0,\ldots, n$.
\end{lemma}
\proof
We deduce from the above discussion that
$$\dist(y^j,x^j)\leq \kappa^n\dist(x,y)\dist(x^j,E)^6\dist(x,E)^{-6}.$$
Since $c_5$ is large, it is enough to show that $n \leq \lof\dist(x,E)Re^{cR}$ for some constant $c>0$. For this purpose, we only have to check that 
$$\dist_P(\xi^j,\xi^{j+1})\geq {e^{-cR}\over \lof\dist(x,E)}\cdot$$ 

By definition of $\xi^j$, we have
$$\epsilon_1\dist(x^j,E)=\dist(x^j,x^{j+1})\leq \dist_P(\xi^j,\xi^{j+1})\max_{t\in[\xi^j,\xi^{j+1}]}\eta(\phi_x(t)).$$
Proposition \ref{prop_eta_E} and Lemma \ref{lemma_close_sing_bis} imply that
$$\eta(\phi_x(t))\simeq  \lof\dist(x^j,E)\dist(x^j,E)\lesssim \lof\dist(x,E)e^{c_3R} \dist(x^j,E).$$
The lemma follows.
\endproof

So, the map $\psi$ is well-defined on $[0,\xi]$. We obtain as in \cite{DinhNguyenSibony2} that it is well-defined on $\overline \D_R$. It is not difficult to see, using the last lemma, that  
$\dist_{\overline\D_R}(\phi_x,\psi)\lesssim e^{-2R}\delta$. 
Since $\phi_y$ is a universal covering map, there is a unique map $\tau:\overline\D_R\to \D$ such that 
$\psi=\phi_y\circ\tau$ and $\tau(0)=0$.
In order to show that $\psi$ satisfies the definition of conformally $(R,\delta)$-close points, 
it remains to check that $\|D\tau\|_\infty\leq A$ and $\|\mu_\tau\|_{\Cc^1}\leq\delta$. It is enough to prove these properties at the point $\xi$ considered above.  

In a neighbourhood of $[\xi^{n-1},\xi^n]$, the maps $\psi$ and $\tau$ are given by
$$\psi=\Phi_n\circ \phi_x \quad\mbox{and}\quad \Phi_n\circ\phi_x=\phi_y\circ\tau.$$
So, we can write locally
$$\tau=\phi_y^{-1}\circ \Phi_n\circ \phi_x.$$

Define $z:=\phi_x(\xi)$ and $w:=\Phi_n(z)$. Then, the norm $\|D\tau(\xi)\|$ with respect to the Poincar\'e metric on $\D$ is bounded by 
$$\eta(z)\|D\Phi_n(z)\|\eta(w)^{-1}.$$
Since $z$ and $w$ are very close to $x^n$, by Proposition \ref{prop_eta_E},
$\eta(z)\eta(w)^{-1}$ is between $1/(2c_1^2)$ and  $2c_1^2$.
Moreover, we have
$$\|\Phi_n-\id\|_{\Cc^2}\lesssim {\dist(\Phi(x^n),x^n)\over \dist(x^n,E)^2}\cdot $$
By Lemma \ref{lemma_delta}, the last quantity is very small. So, we easily deduce that $\|D\tau(\xi)\|\ll 3c_1^2=A$. We also deduce the following useful estimate
$$\|\overline\partial\Phi_n\|_{\Cc^1}\lesssim {\dist(\Phi(x^n),x^n)\over \dist(x^n,E)^2}\cdot $$

It remains to bound $\|\mu_\tau\|_{\Cc^1}$ where $\mu_\tau$ is the Beltrami coefficient defined as in 
\cite{DinhNguyenSibony2}. Recall that this norm is computed with the Euclidean metric on $\D$. 
We first rescale the maps $\phi_x$ and $\phi_y$ using the function $\eta$. This takes into account our distance to the singular points. Define $\zeta:=\tau(\xi)$ and
$$\widetilde\phi_x(t):=\phi_x(\xi+\eta(z)^{-1}t) \quad \mbox{and}\quad 
\widetilde\phi_y(t):=\phi_y(\zeta+\eta(w)^{-1}t).$$
Define also
$$\widetilde\tau:=\widetilde\phi_y^{-1}\circ\Phi_n\circ\widetilde\phi_x.$$

Observe that $\eta(z)$ and $\eta(w)$ are comparable with $\lof \dist(x^n,E)\dist(x^n,E)$ which is bounded from below by $\dist(x^n,E)$. Fix a constant $0<\epsilon_2\ll \epsilon_0$ small enough.
Then, the image of $D_1:=\epsilon_2\dist(x^n,E)\D$ by $\phi_x$ is small and its distance to $E$ is comparable with $\dist(x^n,E)$. Therefore, $\Phi_n(\widetilde\phi_x(D_1))$ is close to $x^n$ and hence it is contained in the image by $\widetilde\phi_y$ of a fixed small disc $D_2$ centered at 0.

The images of $D_1$ and $D_2$ by  $\widetilde\phi_x$ and $\widetilde\phi_y$ are small and contained in some chart of $X$. By Cauchy's formula, we obtain that 
$$\|\widetilde\phi_x\|_{\Cc^2}\lesssim \dist(x^n,E)^{-2}\quad \mbox{and} \quad \|\widetilde\phi_y\|_{\Cc^2}\lesssim \dist(x^n,E)^{-2}.$$ 
Since $\|D\widetilde\phi_y(0)\|=1$ and $\widetilde\phi_y(0)=w$, we deduce that 
$$\|D\widetilde\phi_y^{-1}(w)\|=1 \quad \mbox{and} \quad \|D^2\widetilde\phi_y^{-1}(w)\|\lesssim \dist(x^n,E)^{-2}.$$ 

All these estimates together with the ones on $\|\Phi_n-\id\|_{\Cc^2}$ and $\|\overline\partial\Phi_n\|_{\Cc^1}$ imply that the Beltrami coefficient of $\widetilde\tau$ satisfies the following estimate at 0
$$\|\mu_{\widetilde\tau}\|_{\Cc^1}\lesssim {\dist(\Phi(x^n),x^n)\over \dist(x^n,E)^4}\cdot$$
Recall that $\eta(z)$ and $\eta(w)$ are comparable with $\lof \dist(x^n,E)\dist(x^n,E)$.
Hence, we can deduce from the definitions of $\tau$ and $\widetilde\tau$ the following estimate at the point $\xi$
$$\|\mu_\tau\|_{\Cc^1}\lesssim {\dist(\Phi(x^n),x^n)\over \dist(x^n,E)^6}\cdot$$
By Lemma \ref{lemma_delta}, the last quantity is smaller than $e^{-R}\delta$.
So, the proof of Proposition \ref{prop_conformal_close} is complete for the case where $y$ is the orthogonal projection of $x$ to 
$L_y[3\dist(x,E)\epsilon_0]$.

\bigskip
\noindent
{\bf End of the proof of Proposition \ref{prop_conformal_close}.}
We consider the general case where $y$ is not necessarily equal to 
the orthogonal projection $x'$ of $x$ to 
$L_y[3\dist(x,E)\epsilon_0]$. By hypotheses, $x,y$ are close and satisfy
$$|\log \dist(x,y)|\geq \lof\dist(x,E)e^{c_5R} \quad \mbox{and} \quad 
|\log\dist(x,y)|\geq \lof\dist(y,E)e^{c_5R}.$$ 
We have $\dist(x,x')\leq \dist(x,y)$ and hence $\dist(x',y)\leq 2\dist(x,y)$.

Applying the above construction to $x$ and $x'$, we obtain a map $\psi':\D_R\to L_y$ with $\psi'(0)=x'$ and 
$\dist_{\D_R}(\phi_x,\psi')\ll \delta$ such that the associated map $\tau'$ satisfies 
$$\|D\tau'\|_\infty\ll A \quad \mbox{and}\quad \|\mu_{\tau'}\|_{\Cc^1}\ll\delta.$$
We have to construct a good map $\psi$ associated to $x$ and $y$. 

Observe that $\dist(x',y)\ll \delta\dist(x',E)e^{-e^R}$ since $c_5$ is large. Moreover, $\psi'$ is locally the composition of $\phi_x$ with basic projections from $L_x$ to $L_y$. 
By Proposition \ref{prop_eta_E}, there is a point $a\in \D$ with $|a|\leq \delta e^{-e^R}$ such that $\psi'(a)=y$. So, we can find an automorphism $u:\D_R\to\D_R$ such that $u(0)=a$ and $\|u-\id\|_{\Cc^2}\leq \delta e^{-2R}$. Define $\psi:=\psi'\circ u$. We have $\dist_{\D_R}(\phi_x\circ u,\psi)=\dist_{\D_R}(\phi_x,\psi')\ll\delta$. Since $\eta$ is bounded from above, we also have $\dist_{\D_R}(\phi_x,\phi_x\circ u)\ll\delta$. Therefore, $\dist_{\D_R}(\phi_x,\psi)\ll\delta$. 

The map $\tau$ associated to $\psi$ is given by $\tau:=\tau'\circ u$. It is not difficult to see that $\psi(0)=y$ and 
$$\|D\tau\|_\infty \leq A \quad \mbox{and}\quad \|\mu_\tau\|_{\Cc^1}\leq\delta.$$
So, $\psi$ satisfies the definition of conformally $(R,\delta)$-close points.
\hfill $\square$

\bigskip

We end this section with some technical results that will be used later. Let $m_0$ and $m_1$ be two integers large enough and $\hbar$ be a constant small enough such that $c_1\ll m_0$, $c_1m_0\ll m_1$ and $m_1\hbar\ll \epsilon_0$. 
Define $p:=m_1^4$.

We divide the annulus $\D_{m_1\hbar}\setminus \D_{3\hbar}$ into $40m_1$ equal sectors $S_j$, with $1\leq j\leq 40m_1$, using $40m_1$ half-lines starting at 0 which are equidistributed in $\C$. 
For $x$ in a regular flow box $\U_i$, denote by
$S_l(x)$ the sector of the points $\xi$ in  $\Delta(x,3m_0\hbar)\setminus\Delta(x, m_0\hbar)$ such that $2\pi l/p\leq \arg(\xi)\leq 2\pi (l+1)/p$ for $0\leq l\leq p-1$. Here, $\arg(\xi)$ is defined using the natural coordinate on $\Delta(x,3m_0\hbar)$ centered at $x$.

\begin{lemma} \label{lemma_p_points_1}
Let $x$ be a point in a regular flow box  $\U_i$.  
Then, for every $j$, $\phi_x(S_j)$ contains at least a sector $S_l(x)$. In particular, $\phi_x(S_j)$ intersects
the restriction of $m_1^{-5}\hbar(\Z+i\Z)\times \D^{k-1}$ to a regular flow box $\U_{i'}\simeq \D^k$.
\end{lemma}
\proof
This is a simple consequence of the first assertion in Proposition \ref{prop_eta_E} and that $c_1\ll m_0$ and $c_1m_0\ll m_1$.
Recall that $\phi_x$ is injective on the disc $\D_{\epsilon_0}$ which contains the sectors $S_j$ since $m_1\hbar\ll\epsilon_0$.
\endproof

Consider now a point $x$ in a singular flow box $\U_e$. We deduce from basic properties of the universal covering that there is a unique map $u:\Pi_x\to \D$ such that $\varphi_x=\phi_x\circ u$ and $u(0)=0$.
Fix a constant $t$ satisfying the condition stated just before Lemma \ref{lemma_corr} such that $m_0\hbar<t<2m_0\hbar$.
Define $\zeta_l:= -t\log\|x\|_1 e^{2i\pi l/p}$ and $\xi_l:=u(\zeta_l)$. We use the construction in Section \ref{section_local} for a fixed constant $\lambda$ large enough.
The following two lemmas are related to Lemmas \ref{lemma_corr} and \ref{lemma_J_l_close}.

\begin{lemma} \label{lemma_p_points_2}
Let $x$ be a point in ${3\over 4}\U_e$ satisfying the hypotheses of Lemma \ref{lemma_J_l_close}. Let $w^l$ be as in that lemma. Then, the image $\Gamma$ of the boundary of the disc $D:=-t\log\|x\|_1\D$ by $u$ satisfies $\Gamma\cap \D_{4\hbar}=\varnothing$ and $\Gamma\Subset \D_{m_1\hbar}$. The points $\xi_l$ divide $\Gamma$ into $p$ arcs of length smaller than $m_1^{-2}\hbar$. 
Moreover, for each $l$, there is a point $\xi_l'$ such that $\dist_P(\xi_l',\xi_l)\leq m_1^{-2}\hbar$ and $\phi_x(\xi'_l)=w^l$.
\end{lemma}
\proof
Since $m_0\ll m_1$ and $t\leq 2m_0\hbar \ll\epsilon_0$, by Lemma \ref{lemma_vartheta}, $D$ is contained in the disc of center 0 and of radius $m_1\hbar$ with respect to the Poincar\'e metric on $\Pi_x$. Since holomorphic maps contract the Poincar\'e metric, $u(D)$ is contained in $\D_{m_1\hbar}$.  

The points $\zeta_l$ divide the boundary of $D$ into $p$ arcs.
We deduce from Lemma \ref{lemma_vartheta} that the length of each arc with respect to the Poincar\'e metric, is smaller than 
$4\lambda_*(-\log\|x\|_1)^{-1}$ times its length with respect to the Euclidean metric. 
In particular, it is smaller than $m_1^{-3}\hbar$ since $p=m_1^4$. 
It follows that the $p$ arcs of $\Gamma$ have length smaller than $m_1^{-2}\hbar$. 

Let $\xi$ be a point in the boundary of $\D_{4\hbar}$. The length of $\phi_x([0,\xi])$ with respect to the Poincar\'e metric on $L_x$ is $4\hbar$. By Proposition \ref{prop_eta_E} and Lemma \ref{lemma_vartheta}, 
we have $\eta\lesssim c_1\widehat\eta$ on ${7\over 8}\U_e$. Note that ${7\over 8}\U_e$ contains $\phi_x([0,\xi])$ 
since $\hbar$ is small. Therefore, 
the length of $\phi_x([0,\xi])$ with respect to the Poincar\'e metric on $\widehat L_x=\varphi_x(\Pi_x)$ is bounded by a constant times
$c_1\hbar$. So, the lift of $\phi_x([0,\xi])$ to a curve starting at 0 in $\Pi_x$ is contained in $D$ because $m_0\gg c_1$.  This completes the proof of the first assertion. 

Using again Lemma \ref{lemma_vartheta} and the last assertion in Lemma \ref{lemma_J_l_close}, we obtain that the distance between $z^l$ and $w^l$ with respect to the Poincar\'e metric on $\widehat L_x=\varphi_x(\Pi_x)$ is bounded by a constant times $e^{-23\lambda R}$. By Proposition \ref{prop_eta_E}, this still holds for the Poincar\'e metric on $L_x$. Therefore, there exists a point $\xi_l'$ such that $\dist_P(\xi_l',\xi_l)\leq  m_1^{-2}\hbar$ and $\phi_x(\xi'_l)=w^l$.
\endproof

\begin{lemma} \label{lemma_p_points_3}
Let $x$ be as in Lemma \ref{lemma_p_points_2}. Assume that $x$ is outside the singular flow boxes ${1\over 2}\U_e$. Then, there is a point $\xi_l''$ such that $\dist_P(\xi_l,\xi_l'')\leq m_1^{-4}\hbar$ and $\phi_x(\xi_l'')$ is contained in  the restriction of $m_1^{-5}\hbar (\Z+i\Z)\times \D^{k-1}$ to a
regular flow box $\U_i\simeq \D^k$.
\end{lemma}
\proof
By hypotheses, $x$ belongs to a regular flow box ${1\over 2}\U_i$. Hence, we obtain the result using the same arguments as in the previous lemmas. 
\endproof


\section{Finiteness of entropy: the strategy} \label{section_strategy}

In this section, we will present the strategy for the proof of Theorem \ref{th_main}. 
We have seen in Section \ref{section_local} a simpler situation where we approximate 
Bowen balls by cells. The sizes of these cells, with respect to the Hermitian metric and also with respect to the Poincar\'e metric along the leaves, are very different.  In the global setting, we have considered in Section \ref{section_poincare} the difficulty that a leaf may visit singular flow boxes several times without any obvious rule.

For the proof of Theorem \ref{th_main}, we can imagine that going along the leaves, up to hyperbolic time $R$, transports the cells of the singular flow boxes into the regular part of the foliation and then continue to transport them back to singular flow boxes in different directions. The interaction of the two difficulties increases when $R$ goes to infinity and it is quite hard to get a rough image of the Bowen balls for large time $R$. We are not able to solve this problem, but we can obtain in dimension 2 a good estimate for the number of Bowen balls in the definition of entropy.

We start with a criterium for the finiteness of entropy which allows us 
to reduce the problem to the bound of the entropy of a
finite number of well-chosen transversals.

\begin{proposition} \label{prop_criterium}
Let $\T$ denote the union of the distinguished transversals of the regular flow boxes ${1\over 2}\U_i$. Assume that the entropy of $\T$ is finite. Then, the entropy of $X$ is also finite. 
\end{proposition}

We first prove the following weaker property.

\begin{lemma} \label{lemma_entropy_Y}
Let $Y$ denote the complement of the union of the singular flow boxes ${1\over 2}\U_e$. Then,
under the hypothesis of Proposition \ref{prop_criterium}, the entropy of $Y$ is finite.
\end{lemma}
\proof
Fix a constant $\epsilon>0$. For each flow box $\U_i$, consider a 
$(2R,\epsilon/2)$-dense subset $F_i$ of $\T_i$ with minimal cardinal. Define $F:=\cup F_i$.
Choose a set $\Lambda\subset \D_R$ of less than $e^{10R}$ points which is $e^{-4R}$-dense in $\D_R$, i.e. the discs with centers in $\Lambda$ and with Poincar\'e radius $e^{-4R}$ cover $\D_R$.   
Let $G$ denote the union of the sets $\phi_x(\Lambda)$ with $x\in F$.  
Since $G$ contains at most $e^{10R}\# F$ points, in order to show that the entropy of $Y$ is finite, it suffices to check that any point $z$ in a regular box ${1\over 2}\U_i$ is  $(R,\epsilon)$-close to a point of $G$. 

Consider the plaque $P_z$ of ${1\over 2}\U_i$ which contains $z$. Denote by $y$ the intersection of $P_z$ with the transversal $\T_i$.  So, there is a point $x\in F_i$ such that 
$\dist_{2R}(x,y)\leq \epsilon/2$. Up to a re-parametrization of the leaves, we can assume without loss of generality that $\dist_{2R}(\phi_x,\phi_y)\leq \epsilon/2$. Since $R$ is large, $\phi_y(\D_R)$ contains $P_z$. So, there is a point $\xi\in \D_R$ such that $\phi_y(\xi)=z$. 
It is clear that $z$ and $w:=\phi_x(\xi)$ are $(R,\epsilon/2)$-close. 
Choose a point $\xi'\in \Lambda$ such that $\dist_P(\xi,\xi')\leq e^{-4R}$. 
It is enough to show that $w$ and $w':=\phi_x(\xi')$ are $(R,\epsilon/2)$-close because $w'\in G$.

Observe that there is an automorphism $\tau$ of $\D$ such that $\tau(\xi)=\xi'$ and $\dist_P(\tau(a),a)\leq e^{-R}$ on $\D(\xi,R)$. If $u$ is an automorphism such that $u(0)=\xi$, then $\phi_w:=\phi_x\circ u$ is a covering map of $L_w$ which sends 0 to $w$ and   $\phi_{w'}:=\phi_x\circ \tau\circ u$ is a covering map of $L_{w'}$ which sends 0 to $w'$.
Since the Poincar\'e metric is invariant, we have $\dist_P(\phi_w(a),\phi_{w'}(a))\leq e^{-R}\ll \epsilon$ on $\D_R$. Now, the fact that $\eta$ is bounded from above implies that $w$ and $w'$ are $(R,\epsilon/2)$-close. Hence, the entropy of $Y$ is finite.
\endproof

We will also use the following lemma.

\begin{lemma} \label{lemma_leaf_e}
Let $x$ be a point in $X\setminus E$. Assume that $\phi_x(\D_{2R})$ is contained in a singular flow box ${1\over 2} \U_e$. Let $\epsilon>0$ be a fixed number. If $R$ is large enough, then $\phi_x(\D_R)$ is contained in $\epsilon\U_e$. 
\end{lemma}
\proof
As above, we identify $\U_e$ with $\D^k$. 
Assume that $\phi_x(\D_R)$ contains a point $y$ such that $\|y\|_1\geq \epsilon$. Without loss of generality, assume that the first coordinate $y_1$ of $y$ is a positive number and $\|y\|_1=y_1\geq \epsilon$. By hypothesis, $\phi_y(\D_R)\subset {1\over 2}\D^k$. The real curve $l$ defined by
$$y^t:=\big(t,y_2(t/y_1)^{\lambda_2/\lambda_1},\ldots, y_k(t/y_1)^{\lambda_k/\lambda_1}\big) \quad \mbox{with} \quad t\in [y_1,1/2]$$
is contained in $L_y$ but not in $\phi_y(\D_R)$. Therefore, its Poincar\'e length is at least equal to $R$. 
On the other hand, since $\|y^t\|_1\geq t\geq\epsilon$, we deduce from Proposition \ref{prop_eta_E} that this length is bounded by a constant depending on $\epsilon$. 
This is a contradiction since $R$ is large.
\endproof

\noindent
{\bf End of the proof of Proposition \ref{prop_criterium}.}
Fix a constant $\epsilon>0$ and consider $R>0$ large enough. 
Let $F\subset Y$ be a $(3R,\epsilon/2)$-dense family in $Y$.
Choose a set $W\subset \D_{2R}$ of cardinal $e^{10R}$ which is $e^{-4R}$-dense in $\D_{2R}$. Consider the union $F'$ of the sets $\phi_x(W)$ with $x\in F$. By Lemma \ref{lemma_entropy_Y}, it is enough to show that $F'\cup E$ is $(R,\epsilon)$-dense in $X$, i.e. the Bowen $(R,\epsilon)$-balls centered at a point in $F'$, cover $X$. 

Consider a point $z\in X$. If $\phi_z(\D_{2R})$ is contained in a singular flow box ${1\over 2}\U_e$, by Lemma \ref{lemma_leaf_e}, we have 
$\dist_R(e,z)< \epsilon$. So, for $R$ large enough, $z$ belongs to the Bowen $(R,\epsilon)$-ball centered at $e$.  It remains to consider the case where $\phi_z(\D_{2R})$ contains a point $y\in Y$. 

Let $x$ be a point in $F$ such that $\dist_{3R}(x,y)\leq \epsilon/2$. So, we can find covering maps $\phi_x$ and $\phi_y$ such that $\dist_{3R}(\phi_x,\phi_y)\leq \epsilon/2$. Since $y\in\phi_z(\D_{2R})$, there is a point $\xi\in \D_{2R}$ such that $\phi_y(\xi)=z$. It is clear that $z$ and $w:=\phi_x(\xi)$ are $(R,\epsilon/2)$-close. Let $\xi'$ be a point in $W$ such that $\dist_P(\xi,\xi')\leq e^{-4R}$. The point $w':=\phi_x(\xi')$ belongs to $F'$. We show as in Lemma \ref{lemma_entropy_Y} that $w$ and $w'$ are $(R,\epsilon/2)$-close. Therefore, $z$ belongs to the Bowen $(R,\epsilon)$-ball of center $w'$.
\hfill $\square$

\medskip

From now on, assume that $\dim X=2$ and hence $\dim \T=1$. 
Our strategy  is to construct an adapted covering of $\T$ by open discs. 
Consider a hyperbolic time $R$ large enough. For simplicity, assume that $R=N\hbar$ with $N$ integer.  We will construct a set $\widetilde \T$ which contains $\T$ and other transversals ouside and inside the singular flow boxes. Then, we will construct by induction on $m$, with $2m_1\leq m\leq N$, a covering $\Vc_m$ of $\widetilde\T$ by discs such that if two points $x, y$ belong to the same disc, the associated leaves $L_x$ and $L_y$ are close until time $m\hbar$. Of course, we have to control the cardinal of $\Vc_m$ in order to deduce the finiteness of entropy of $\T$ by taking $m=N$.

The covering $\Vc_m$ is obtained by refining a finite number of other coverings of $\widetilde \T$. We will use the following technical lemma in order to estimate the cardinal of $\Vc_m$ and to show that the cardinal of $\Vc_N$ grows at most exponentially when $N$ tends to infinity.

\begin{lemma}\label{lemma_disc_count}
Let $K$ be a finite family of sets such that each of them is contained in a complex plane, i.e. a copy of $\C$. 
Let $\Vc^i$ with $1\leq i\leq n$ be $n$ coverings of $K$ by  less than $M$ discs.
Then, we can cover  $K$ with a family $\Vc$ of less than
$200^nM$ discs such that $\Vc\prec \Vc^i$ in the sense that
every disc $D\in\Vc$ satisfies $2D\subset   2D_1\cap \cdots\cap 2D_n$ for some
$D_i\in \Vc^i$.
 \end{lemma}
\proof
By induction, it is enough to consider the case $n=2$ and to find a covering $\Vc$ with less than $200M$ discs.
The covering $\Vc$ contains two kinds of discs that we construct below.

An element $D_2$ of $\Vc^2$ belongs to $\Vc$ iff there is  $D_1\in\Vc^1$ such that $D_1\cap D_2\not=\varnothing$ and  $\radius(D_1)>2\radius(D_2)$.
Clearly, these discs satisfy the last condition in the lemma. 

Consider now $D_j\in\Vc^j$ such that $D_1\cap D_2\not=\varnothing$ and $\radius(D_1)\leq 2\radius(D_2)$. Denote by $4\rho_1$ the radius of $D_1$.
A disc $D$ of the second kind is a disc of radius $\rho_1$ centered at a point in $\rho_1(\Z+i\Z)$ which intersects $D_1\cap D_2$. 
It is clear that $2D\subset 2D_1\cap 2D_2$. Note that such a disc $D$ can be associated to several discs $D_2$. 

Observe that each disc $D_1$ as above is associated to less than 100 discs of the second kind. Therefore,
$\Vc$ contains less than $200M$ discs. This family covers $K$ since it covers $D_1\cap D_2$ for all $D_j\in \Vc^j$.  This completes the proof.
 \endproof

As we mentioned above, the covering $\Vc_m$ will be obtained by induction on $m$. 
We will construct $p$ coverings $\Vc_m(l)$ of $\widetilde\T$ by discs which are obtained using the images of discs in $\Vc_{m-1}$ by some holonomy maps (the integer $p=m_1^4$ was fixed above). Near a leaf, such a holonomy map looks like a displacement following a given direction with a small hyperbolic time.  The covering $\Vc_m$ will be obtained by refining the $p+1$ coverings $\Vc_{m-1}$ and $\Vc_m(l)$, with $0\leq l\leq p-1$, thanks to Lemma \ref{lemma_disc_count}. The details will be given in Section \ref{section_transversal}.

A crucial property of $\Vc_m$ is that when two points $x,y$ belong to the same element of $\Vc_m$, the maps $\phi_x$ and $\phi_y$ are close on $\D_{m\hbar}$. This property will be obtained by induction, i.e. from a similar property of $\Vc_{m-1}$. For this purpose, we need to cover $\D_{m\hbar}$ and $\phi_y(\D_{m\hbar})$ by discs of radius $(m-1)\hbar$ in order to apply the induction argument. 
The following lemma shows that we only need a fixed number of such discs.

\begin{lemma} \label{lemma_cover_D}
Let $\xi_j$ be a point in $S_j$ with $1\leq j\leq 40m_1$. Then, for $m\geq  2m_1$, the disc $\D(0,(m-1)\hbar)$ and the $40m_1$ discs $\D(\xi_j,(m-2)\hbar)$ cover the disc $\D_{m\hbar}$.  
\end{lemma}
\proof
Let $\xi$ be a real number such that $3\hbar\leq \xi<2m_1\hbar$. We claim that it is sufficient to show that $\D(\xi,(m-2)\hbar)$ contains all points $\zeta\in \D_{m\hbar}\setminus\D_{(m-1)\hbar}$ such that $|\arg(\zeta)|\leq \pi/(20m_1)$.  Indeed, this property implies that $\D(\xi_j,(m-2)\hbar)$ contains the sector 
$|\arg(\zeta)-\arg(\xi_j)|\leq \pi/(5m_1)$ in $\D_{m\hbar}\setminus\D_{(m-1)\hbar}$. The union of these sectors covers $\D_{m\hbar}\setminus\D_{(m-1)\hbar}$. 

Denote by $\zeta'$ the intersection of the half-line through $\zeta$ started at 0 with the circle of center 0 through $\xi$. Since $2m_1\hbar$ is small, it is not difficult to see that $\dist_P(\xi,\zeta')<\hbar$. Moreover, we have
$$\dist_P(\zeta,\zeta')=\dist_P(0,\zeta)-\dist_P(0,\xi)\leq (m-3)\hbar$$ 
Therefore, 
$$\dist_P(\xi,\zeta)\leq \dist_P(\xi,\zeta')+\dist_P(\zeta,\zeta')<(m-2)\hbar.$$ 
The lemma follows.
\endproof

We deduce from the last lemma the following result which is more adapted to our problem when we are near a singular point.

\begin{lemma} \label{lemma_cover_D_bis}
Let $\Gamma$ be a closed curve in $\D_{m_1\hbar}\setminus \D_{4\hbar}$ such that $0$ does not belong to the unbounded component of $\C\setminus\Gamma$. We divide it into $p$ arcs $\wideparen{\xi_j\xi_{j+1}}$ with $0\leq j\leq p$ and $\xi_p=\xi_0$. Assume that the length of 
$\wideparen{\xi_j\xi_{j+1}}$ with respect to the Poincar\'e metric on $\D$ is smaller than $m_1^{-2}\hbar$. 
Let $\xi_j'$ be a point in $\D$ such that $\dist_P(\xi_j,\xi_j')\leq m_1^{-2}\hbar$.
Then, the  disc
$\D(0,(m-1)\hbar)$ and the $p$ discs $\D(\xi'_j,(m-1)\hbar)$ cover the disc $\D_{m\hbar}$.  
\end{lemma}
\proof
Observe that $\Gamma\cap S_j$ admits a point $\xi$ such that the disc $\D(\xi,3m_1^{-2}\hbar)$ is contained in $S_j$. If $\xi$ belongs to $\wideparen{\xi_i\xi_{i+1}}$, then this arc and also the point $\xi'_i$ are contained in $S_j$. Lemma \ref{lemma_cover_D} implies the result.
\endproof


\section{Adapted transversals and their coverings} \label{section_transversal}

In this section, assume that $X$ is a compact complex surface. We will construct, for every integer $N$ large enough, the family of transversals $\widetilde\T=\T^\reg\cup \T^\sing$ for $R=N\hbar$ and its coverings $\Vc_m=\Vc_m^\reg\cup \Vc_m^\sing$ with $2m_1\leq m\leq N$. We will use the constants introduced in the list  given in the introduction.

Recall that each regular flow box $\U_i$ is identified to $\D^2$ and is associated with the distinguished transversal $\{0\}\times \D$. Define $\T_i^a:=\{a\}\times \D$ with $a$ in the lattice $m_1^{-5}\hbar(\Z+i\Z)\cap\D$.  Denote by $\T^\reg$ the family of these transversals. Note that this family does not depend on $R$.
We will, however, consider each element $\T_i^a$ of $\T^\reg$ with multiplicity $e^{46\lambda R}$ and 
we denote them by $\T_i^a(1),\ldots, \T_i^a(e^{46\lambda R})$. These transversals are considered as distinct. 
We will construct later the covering $\Vc_m$ of $\widetilde \T$ by induction on $m$ using 
Lemma \ref{lemma_disc_count} applied to 
$\Vc_{m-1}$ and $p$ other coverings of $\widetilde\T$. 
The discs used to cover $\T_i^a(j)$  depend on the index $j$. 
The multiplicities allow us to get a good bound for the number of discs in $\Vc_m$, see Proposition \ref{prop_cover_count} below.   

Consider now a singular flow box $\U_e\simeq \D^2$ and define $\T^a_e:=\{a\}\times{3\over 4}\D$ with $a\in\Sigma\cap {3\over 4}\D$. Recall that $\Sigma$ was constructed in Section \ref{section_local} and we use here a large fixed constant $\lambda$. Denote by $\T^\sing$ the family of these transversals $\T^a_e$, where each element is counted only one time. Note that $\T^\sing$ depends on $R$. Define $\widetilde\T$ the union of $\T^\reg$ and $\T^\sing$.

We now construct the covering $\Vc_{2m_1}$ of $\widetilde\T$. 
Choose for $\T^\reg$ a covering $\Vc^\reg_{2m_1}$ by less than 
$e^{70\lambda R}$ discs of radius  $e^{-10R}$, where we count the multiplicities of transversals.
Consider the family $\Vc^\sing_{2m_1}$ of the discs in $\T_e^a$ centered at $(a,b)\in\Sigma^2$ with radius $100e^{-23\lambda R}|b|$ if $b\not =0$ and of radius $\alpha_1$ if $b=0$. It is not difficult to check that this family covers $\T^\sing$. Define $\Vc_{2m_1}$ as the union of $\Vc_{2m_1}^\reg$ and $\Vc_{2m_1}^\sing$. The total number of discs used here is bounded by $e^{200\lambda R}$. Recall that we only consider a large $R$. So, we have $m_1\ll R$.

For each point $x\in \T_i^a$ and $\rho$ small enough, denote by $D(x,\rho)$ the disc of center $x$ and of radius $\rho$ in $2\T_i^a$ and $\Delta(x,\rho)$ the disc of center $x$ and of radius $\rho$ in the plaque of $2\U_i$ containing $x$.  The notation $D(x,\rho)$ can be used for $x$ in a transversal $\T_e^a$.
Since $m_1\hbar\ll \epsilon_0$, we have the following useful
properties for large $R$, where we use that the discs are of size less than $e^{-10R}$:

\begin{enumerate}
\item[(H1)] If $x,y$ are in $2D$ for some $D$ in $\Vc^\reg_{2m_1}$, then the basic projection $\Phi$ associated to $x$ and $y$ exists and sends $\Delta(x,2m_1\hbar)$ to $\Delta(y,6m_1\hbar)$. Moreover, we have $\|\Phi-\id\|_{\Cc^2}\leq e^{-9R}$ on $\Delta(x,2m_1\hbar)$. 
\item[(H1)'] If $x,y$ are in $2D$ for some $D$ in $\Vc^\sing_{2m_1}$, then the basic projection $\Phi$ associated to $x$ and $y$ exists as in Lemma \ref{lemma_proj_E}. It  satisfies $\|\Phi-\id\|_{\Cc^0}\leq e^{-9R}$ on $L_x\cap {3\over 4}\D^2$. 
\item[(H2)] Consider two transversals $\T_i^a$ and $\T_j^b$.
If $x$ is in $\T_i^a$ such that $\Delta(x,2m_1\hbar)$ intersects $\T_j^b$ at a point $y$, then the holonomy map $\pi$ from $2\T_j^b$ to $2\T_i^a$ is well-defined on $D(y,4\hbar)$  with image in ${3\over 2}\T_i^a$. 
\item[(H3)] If $D$ is a disc contained in $D(y,\hbar)$, then $\pi(D)$ is {\it quasi-round}, i.e. there is a disc $D'\subset {3\over 2}\T_i^a$ such that $D'\subset \pi(D)\subset {11\over 10}D'$ and $2D'\subset \pi(2D)$.
\end{enumerate}
For the property (H3), we use the fact that the holonomy map is holomorphic with no critical point. So, on small discs,  it is close to homotheties.

The construction of $\Vc_m=\Vc_m^\reg\cup \Vc_m^\sing$ will be obtained by induction on $m$. 
It only contains small discs of diameter less than $e^{-10R}$. 
Assume that the construction is done for $m-1$. 
In order to obtain $\Vc_m$, we will apply Lemma \ref{lemma_disc_count} to $K:=\widetilde \T$, to the covering $\Vc_{m-1}$ and $p$ other coverings $\Vc_m(l)=\Vc_m^\reg(l)\cup \Vc_m^\sing(l)$ with $0\leq l\leq p-1$.
Lemma \ref{lemma_disc_count} allows us to obtain a covering $\Vc_m$ such that 
$\Vc_m\prec \Vc_{m-1}$ and $\Vc_m\prec \Vc_m(l)$. Roughly speaking, we will cover each disc in $\Vc_{m-1}$ by smaller discs, so that after traveling a fixed time in some direction starting from such a small disc, we arrive at a disc of $\Vc_{m-1}$ in another transversal.

We explain now the construction of  $\Vc_m^\reg(l)$.  There are two cases to consider. Recall that $S_l(x)$ is defined just before Lemma \ref{lemma_p_points_1}.

\medskip
\noindent
{\bf Case 1a.}
Assume that there is a point $x_0\in \T_i^a$ such that $S_l(x_0)$ intersects a singular flow box ${1\over 2}\U_e$.  Consider an arbitrary point $x\in\T_i^a$. 
Since we are still 
far from singular points, $S_l(x)$ contains and is contained in small discs of size independent of $R$. 
Thus, there are about a constant times $e^{46\lambda R}$ transversals $\T_e^b$ which intersect $S_l(x)$.
These transversals $\T_e^b$ are still far from the singularities.  So, the following properties hold because the regular flow boxes are of small size and $\hbar$ is small:

\begin{enumerate}
\item[(H2)'] There is a well-defined holonomy map from $2\T_i^a$  to $\T_e^b$. Denote by $\pi$ its inverse.
\item[(H3)'] If $D$ is a disc in $\T_e^b$ of radius less than $\hbar$ which intersects $\pi^{-1}(\T_i^a)$, then $\pi$ is defined on $2D$ and $\pi(D)$ is quasi-round, i.e. there is a disc $D'\subset 2\T_i^a$ such that $D'\subset \pi(D)\subset {11\over 10}D'$ and $2D'\subset \pi(2D)$. 
\end{enumerate}

The property (H3)' allows us to cover $\pi(D)$ by $D'$ and its $100$ {\it satellites} which are 100 discs $D_n'$, $0\leq n\leq 99$, 10 times smaller than $D'$, and such that $D'_n\cap D'\not=\varnothing$. Notice that $2D_n'\subset \pi(2D)$ for all $n$. We will use this important property later.
So, we can cover each $\T_i^a(s)$ with the discs $D'$ obtained above for $D\in\Vc_{m-1}$ together with theirs satellites $D_n'$. The choice of $\T_e^b$ depends on $s$. Thanks to the multiplicities of the $\T_i^a$, we can use each $\T_e^b$ only a bounded number of times.
The reason for introducing those multiplicities is that in the intersection of regular and singular flow boxes the $\T_e^b$ are more dense than the $\T_i^a$. 

\medskip
\noindent
{\bf Case 1b.} Assume that $S_l(x)$ does not intersect any singular 
flow box ${1\over 2}\U_e$ for every $x\in\T_i^a$. 
For all holonomy maps $\pi$ and discs $D$ in $\Vc_{m-1}$ satisfying (H2) and (H3) for some $x$ and $y$ such that $y\in S_l(x)$, we choose a disc $D'$ as in (H3) and 100 satellites of $D'$ in order to cover $\pi(D)$. Note that for any $x\in \T_i^a$, there exists a choice of $y,\pi,D$ such that $x\in \pi(D)$, see also Lemma \ref{lemma_p_points_1}. It follows that the construction gives us a covering of $\T_i^a(s)$ using the elements of $\Vc_{m-1}$ which cover the transversals $\T_i^b(s)$. We make sure to use here the same index $s$, in particular, each disc in $\Vc_{m-1}$ is used a bounded number of times. 
This ends the construction of the covering $\Vc_m^\reg(l)$ of $\T^\reg$. 

\bigskip

The construction of $\Vc_m^\sing(l)$ is more delicate. First, we always add to $\Vc_m^\sing(l)$ the disc of center $(a,0)$ and of radius $\alpha_1$ in $\T_e^a$ and call it {\it an exceptional disc}.
If $|a|\leq\alpha_1$, we just choose $\Vc_m^\sing(l)$ equal to $\Vc_{m-1}$ on $\T_e^a$. 
In this case, we also say that $\T_e^a$ is an {\it  an exceptional transversal}.
Consider now a transversal $\T_e^a$ with $\alpha_1<|a|\leq 3/4$. 
Consider a point $x=(a,d)\in \T_e^a\cap \Sigma^2$ with $|d|\geq \alpha_1$. 
Recall that the map $J_l$ is defined just before Lemma \ref{lemma_corr}.
By Lemma \ref{lemma_close_sing_bis}, the point $w:=J_l(x)$ satisfies the hypotheses of Lemmas \ref{lemma_J_l_close} and \ref{lemma_p_points_2} since $\lambda$ is a large constant.
We also distinguish two cases. 

\medskip\noindent
{\bf Case 2a.} Assume that $w$ belongs to ${5\over 8}\U_e$ (note that ${5\over 8}\U_e\subset {3\over 4}\U_e$). Consider the transversal $\T_e^b$ which contains $w$ and write $w=(b,v)$. If $D\subset \T_e^b$ is an element of $\Vc_{m-1}^\sing$ such that $\dist (w,D)\leq 100\gamma |v|e^{-23\lambda R}$, denote by $D'$ the disc on $\T_e^a$ which is the image of $D$ by $\Psi_{w,x}$ since this map preserves $\Sigma^2$. By Lemma \ref{lemma_J_l_close},  the obtained discs $D'$ cover the disc of center $x$ and of radius $100|d|e^{-23\lambda R}$ in $\T_e^a$. They are elements of $\Vc_m^\sing(l)$. Note that we use here the property that $\Psi_{w,x}$ is conformal. This is the only point where the hypothesis $\dim X=2$ is essential. 

\medskip\noindent
{\bf Case 2b.} Assume that $w$ is not in ${5\over 8}\U_e$. Since $\hbar$ is small, by Lemma \ref{lemma_p_points_2}, $x, w$ are outside ${1\over 2}\U_e$ and $\varphi_x(-t\log\|x\|_1\D)$ is contained in a plaque of a regular flow box. 
By Lemma \ref{lemma_p_points_3}, we can find a transversal $\T_i^b$ which intersects $\varphi_x(-t\log\|x\|_1\D)$ at a point $y$ near $\varphi_x(-t\log\|x\|_1e^{2i\pi l/p})$, i.e. the distance between these two points is less than $m_1^{-3}\hbar$. We have the following properties: 
\begin{enumerate}
\item[(H2)''] There is a well-defined holonomy map $\pi$ from $2\T_i^b$  to $\T_e^a$. 
\item[(H3)''] If $D$ is a disc in $2\T_i^b$ of radius less than $\hbar$ which intersects $\T_i^b$, then $\pi$ is defined on $2D$ and $\pi(D)$ is quasi-round in the sense that there is a disc $D'\subset \T_e^a$ such that $D'\subset \pi(D)\subset {11\over 10}D'$ and $2D'\subset \pi(2D)$. 
\end{enumerate}

We cover a neighbourhood of $x$ with the discs $D'$ and theirs satellites as above with $D\in\Vc_{m-1}$. 
We make sure that for each $\T_e^a$ we only use discs from $\T_i^b(s)$ for a fixed $s$. 
It is important to observe that the last construction concerns about a constant times $e^{46\lambda R}$ transversals $\T_e^a$. Therefore, we can choose the index $s$ so that 
each disc in $\Vc_{m-1}^\reg$ is used a bounded number of times.  
We also fix a choice which does not depend on $m$.
This ends the construction of the covering $\Vc_m^\sing(l)$. 

\medskip

We have the following crucial proposition.

\begin{proposition}\label{prop_cover_count}
There is a constant $c>1$ independent of $R$ such that the cardinal of $\Vc_m$ is smaller than $c^R$ for 
$2m_1\leq m \leq N$. 
\end{proposition}
\proof
By Lemma \ref{lemma_disc_count} applied to $n:=p+1$, it is enough to show that in the above construction of $\Vc_m(l)$, each disc in $\Vc_{m-1}$ is used less than $c'$ times where $c'$ is a constant. This can  be checked step by step in our construction. For one of these steps, we use Lemma \ref{lemma_corr}. We obtain by induction that
$$\# \Vc_m \leq (200^{p+1}c')^{m-2m_1} \# \Vc_{2m_1} \leq (200^{p+1}c')^{m-2m_1} e^{200\lambda R}.$$
This implies the proposition.
\endproof

The last proposition shows that the cardinal of $\Vc_N$ is smaller than $c^R=c^{N\hbar}$.
Consider a disc $D$ in $\Vc_N$. It is constructed by induction  
using the holonomy maps $\pi$ as in (H2), (H2)', (H2)'' or the map $\Psi_{w,x}$ as above. This corresponds to 
the 4 cases described above. 
By Proposition \ref{prop_criterium}, in order to obtain Theorem \ref{th_main},  we will show in Proposition \ref{prop_end} below that two points in the same disc $D$ are $(R/3, e^{-R/4})$-close. This will be done using the notion of conformally $(R,\delta)$-close points. 

We will associate to a disc $D$ as above a tree $F_D$ which partially encodes the construction of $\Vc_N$. 
Its combinatorial and metric properties (see Lemma \ref{lemma_discrete_leaf} below) will allow us  to construct conformally $(R,\delta)$-close maps from leaves to leaves following the tree.
Points in $D$ are associated to some isomorphic trees and the isomorphisms are coherent with the dynamics of the foliation. 
 
The set of vertices of the tree $F_D$ 
will be the union $F_D(0)\cup F_D(1)\cup\ldots\cup F_D(N-2m_1)$, where $F_D(0)=\{D\}$ and $F_D(m)\subset \Vc_{N-m}$. Moreover, each point in $F_D(m)$ is joined to a unique point in $F_D(m-1)$ and to at most $p$ points in $F_D(m+1)$. We give now the construction of $F_D$ by induction.

If $D$ belongs to an exceptional disc or an exceptional transversal, we just take $F_D(1)=\cdots=F_D(N-2m_1)=\varnothing$. Otherwise, $D$ is obtained using $p$ holonomy maps $\pi_i$ as in (H2), (H2)', (H2)''  or the map $\Psi_{w,x}$ as above. By construction, there are $D_i$ in $\Vc_{N-1}$ such that $2D$ is contained in $\pi_i(2D_i)$. We choose $F_D(1):=\{D_1,\ldots,D_p\}$. Each $D_i$ is joined to $D$. We then obtain a part of the tree.

In order to obtain $F_D(2)$, we will repeat the above construction but for each $D_i\in F_D(1)$ instead of $D$. If $D_i$ belongs to an exceptional disc or an exceptional transversal, then it is not joined to any element in $F_D(2)$. Otherwise, it is joined to $p$ elements in $F_D(2)\subset \Vc_{N-2}$. 
We then continue the same construction in order to obtain $F_D(3),\ldots,F_D(N-2m_1)$. Note that $F_D$ is not uniquely determined by $D$ but we fix here a choice for each disc $D$. 

Now, we construct for each point $x\in 2D$ a tree $F_x\subset \D$ which is canonically isomorphic to $F_D$. 
The set of vertices of $F_x$ will be $F_x(0)\cup F_x(1)\cup\ldots\cup F_x(N-2m_1)$ such that $F_x(0)=\{0\}$ and $\phi_x$ sends each point in $F_x(m)$ to a disc $2D'$ with $D'\in F_D(m)$ and defines a bijection between $F_x(m)$ and $F_D(m)$. Moreover, a point in $F_x(m+1)$ and a point in $F_x(m)$ are joined by an edge if and only if the associated vertices in $F_D$ are also joined by an edge. We have to give here some details because $\phi_x$ is not injective in general.

We obtain $F_x(1)$ as follows. If $D$ belongs to an exceptional disc or an exceptional transversal, then we take 
$F_x(1)=\varnothing$. Otherwise, we distinguish two cases.
When $D$ is outside the singular flow boxes ${1\over 4}\U_e$, then $\phi_x$ is injective on $\D_{\epsilon_0}$ and the image of this disc intersects $2D_i$ at a unique point for any $D_i\in F_D(1)$. Therefore, it is enough to define $F_x(1)$ as the pull-back by $\phi_x$ of these intersection points in $\D_{\epsilon_0}$. 
Consider now the case where $D$ intersects a singular flow box ${1\over 4}\U_e$. By construction, there is a set $G$ of  $p$ points very close to the points $-t\log\|x\|_1e^{2i\pi l/p}$ in $\Pi_x$ which are sent by $\varphi_x$ to the discs in $F_D(1)$. If $\tau:\Pi_x\to \D$ is the unique holomorphic map 
such that $\tau(0)=0$ and $\varphi_x=\phi_x\circ \tau$, define $F_x(1):=\tau(G)$. Clearly, $\phi_x$ defines a bijection between $F_x(1)$ and $F_D(1)$.  Each point in $F_x(1)$ is joined to $F_x(0)=\{0\}$. 

In order to obtain $F_x(2)$, it is enough to repeat the same construction to each point $\phi_x(a)$ with $a\in F_x(1)$. We use that $\phi_x$ is injective on $\D(a,\epsilon_0)$ if $\phi_x(a)$ is outside the singular flow boxes ${1\over 4}\U_e$ and otherwise there is a unique holomorphic map $\tau:\Pi_x\to \D$
such that $\tau(0)=a$ and $\varphi_a=\phi_x\circ \tau$. By induction, we obtain $F_x(3),\ldots, F_x(N-2m_1)$ satisfying 
the properties stated above.

The following lemma is essential for the proof of Theorem \ref{th_main}.

\begin{lemma}  \label{lemma_discrete_leaf}
The set of vertices $F_x(m)$ is contained in $\D_{m_1m\hbar}$ for every $m$. If
$a$ is a point in $F_x$ which is joined to $p$ points $a_1,\ldots, a_p$, then the union of $\D(a,(m-1)\hbar)$ and 
$\D(a_i,(m-1)\hbar)$ contains $\D(a,m\hbar)$ for every $m\geq 2m_1$. 
\end{lemma}
\proof
We prove the first assertion. It is enough to check that if a point $\zeta$ in $F_x(n-1)$ is joined to a point $\xi$ in $F_x(n)$, then $\dist_P(\xi,\zeta)\leq m_1\hbar$. For simplicity, we consider the case where $n=1$, the general case is obtained in the same way. So, we have $\zeta=0$ and $\xi$ is in $F_x(1)$. 
Define $y:=\phi_x(\xi)$. This is the preimage of $x$ by a holonomy map $\pi$ as in (H2), (H2)', (H2)'' or by a map $\Psi_{w,x}$ as above. 

Now, if $y$ is given by $\pi$ as in (H2),   the distance between $x$ and $y$ is smaller than $3m_0\hbar$. By Proposition \ref{prop_eta_E}, the Poincar\'e distance between 0 and $\xi$ is smaller than $m_1\hbar$ since $m_1\gg c_1m_0$.

For the other cases, $y$ is very close to a point 
$y':=\varphi_x(-t\log\|x\|_1e^{2i\pi l/p})$ with respect to the Poincar\'e metric. Recall that $m_0\hbar<t<2m_0\hbar\ll 1$. The Poincar\'e distance between $x$ and $y'$ 
is bounded by the distance between 0 and $-t\log\|x\|_1 e^{2i\pi l/p}$ with respect to the Poincar\'e metric on $\Pi_x$. Therefore, by Lemma \ref{lemma_vartheta}, we deduce that $\dist_P(x,y')$ is bounded by a constant times $m_0\hbar$. It follows easily that
 $\dist_P(0,\xi)$, which is equal to $\dist_P(x,y)$, is smaller than $m_1\hbar$. This completes the proof of the first assertion.

We prove now the second assertion.  By Lemmas \ref{lemma_cover_D} and \ref{lemma_cover_D_bis}, we only have to check that $\{a,a_1,\ldots,a_p\}$ contains a subset satisfying the hypotheses of those lemmas.
If $a$ is in $\T^\reg$, this is a consequence of Lemma \ref{lemma_p_points_1}. If $a$ is in $\T^\sing$,  Lemmas \ref{lemma_p_points_2} and \ref{lemma_p_points_3} imply the result. 
\endproof


\section{Finiteness of entropy: end of the proof} \label{section_proof}

In this section, we complete the proof of Theorem \ref{th_main}. By Proposition \ref{prop_criterium}, it is enough to show that the entropy of $\T^\reg$ is finite. By Proposition \ref{prop_cover_count}, we only have to check that each disc in $\Vc_N^\reg$ is contained in a Bowen $(R/3,e^{-R/4})$-ball. So, Theorem \ref{th_main} is a consequence of the following proposition.

\begin{proposition} \label{prop_end}
Let $x,y$ be two points of a disc $2D\subset 2\T_i^a$ with $D$ in $\Vc_N$. Then, they are $(R/3,e^{-R/4})$-close.
\end{proposition}

The proof of this result uses Proposition \ref{prop_quasi_close} and occupies the rest of this section. Recall that by Proposition \ref{prop_eta_E}, we have $c_1^{-1}\leq |\eta|\leq c_1$ on $\T^\reg$. 
So, it is enough to check that $x, y$ are conformally $(R,e^{-3R})$-close and we have to construct a map $\psi$ from $\D_R$ to $L_y$ which is close to $\phi_x$ as in the definition of conformally $(R,e^{-3R})$-close points. Proposition \ref{prop_end} is a direct consequence of Lemmas \ref{lemma_end_1} and \ref{lemma_end_2} below, which correspond respectively to the case where the leaves are far from the separatrices and to the case where the leaves are close to some separatrice.

We have the following lemma which holds for $x,y$ in $2D$ with an arbitrary disc $D$ in $\Vc_N=\Vc_N^\reg\cup\Vc_N^\sing$. 
Recall that the trees $F_x$ and $F_y$ are isomorphic to $F_D$. Therefore, there is an isomorphism $\sigma$ from $F_x$ to $F_y$.
Denote by $x'$ the image of $x$ by the basic projection associated to $x$ and $y$. 
 We say that the tree $F_D$ is {\it complete} if each element in $F_D(m-1)$ is joined to $p$ elements in $F_D(m)$ for any $1\leq m\leq N-2m_1$.
 
\begin{lemma} \label{lemma_end}
Assume that the tree $F_D$ is complete. Then, there is a unique map $\psi':\D_{N\hbar}\to L_y$ which is locally 
 the composition of $\phi_x$ with basic projections from leaves to leaves and is equal, in a neighbourhood of any point $a$ in $F_x\cap \D_{N\hbar}$, to the composition of $\phi_x$ with the basic projection associated to $\phi_x(a)$ and $\phi_y(\sigma(a))$. 
\end{lemma}
\proof
Observe that by continuity, if such a map $\psi'$ exists, it is unique.
We show by induction on $m$ that such a projection exists on $\D(a,{m\hbar})$ for $a$ in $F_x(0)\cup\ldots \cup F_x(N-m)$ with $2m_1\leq m\leq N$. Define $\widetilde x:=\phi_x(a)$ and $\widetilde y:=\phi_y(\sigma(a))$. Denote also by $\widetilde D$ the element of $F_D$ associated to $\widetilde x$ and $\widetilde y$. This is an element of 
$F_D(0)\cup\ldots \cup F_D(N-m)$ such that $\widetilde x$ and $\widetilde y$ belong to $2\widetilde D$.

If $\widetilde x$ or $\widetilde y$ is outside the singular flow boxes ${1\over 2}\U_e$, since $2m_1\hbar$ is small, $\phi_x$ sends $\D(a,{2m_1\hbar})$ bijectively to a disc in $L_{\widetilde x}[\epsilon_0]$. So, it is not difficult to see that the desired property holds in this case for $m=2m_1$. When $\widetilde x$ and $\widetilde y$ belong to a singular flow box ${1\over 2}\U_e$, the property for $m=2m_1$ is a consequence of Lemma \ref{lemma_proj_E}. 
Assuming now the property for $m-1$, we have to show it for $m$. 

Recall that $2\widetilde D$ is contained in $2D'\cap \pi_0(2D_0)\cap \ldots\cap \pi_{p-1}(2D_{p-1})$. Here, $D',D_i$ are elements of $\Vc_{m-1}\cup \ldots \cup \Vc_{N-1}$ and the $\pi_i$ are holonomy maps from some transversals in $\T$ to the one containing $\widetilde D$ or a map $\Psi_{w,x}$ as in Section \ref{section_transversal}. Define $x^i:=\pi_i^{-1}(\widetilde x)$ and $y^i:=\pi_i^{-1}(\widetilde y)$. Denote by $\xi^i$ the points in $F_x$ such that $\phi_x(\xi^i)=x^i$. 

It follows from the induction hypothesis  that there is a map $\psi_i'$ from $\D(\xi^i,(m-1)\hbar)$ to $L_y$ which is the composition of $\phi_x$ with basic projections from leaves to leaves. Moreover, it is equal, in a neighbourhood of any point $b$ in $F_x\cap  \D(\xi^i,(m-1)\hbar)$, to the composition of $\phi_x$ with the basic projection associated to $\phi_x(b)$ and $\phi_y(\sigma(b))$.  Also by induction hypothesis, there is an analogous map defined on 
$\D(a,{(m-1)\hbar})$. From the uniqueness of these maps, we can glue them together. 
Using the second assertion of Lemma \ref{lemma_discrete_leaf}, we obtain a map $\psi'$ defined on $\D(a,{m\hbar})$.
\endproof

\begin{lemma} \label{lemma_end_1}
Let $x,y,D$ be as in Proposition \ref{prop_end}. Assume that the tree $F_D$ is complete. Then, $x$ and $y$ are conformally $(R,e^{-3R})$-close.
\end{lemma}
\proof
We have to construct a map $\psi$ satisfying the definition of conformally $(R,e^{-3R})$-close points.
Let $\psi'$ be the map constructed in  Lemma \ref{lemma_end}. We deduce from the construction using (H1), (H1)'  and Lemma \ref{lemma_proj_E} that $\|\psi'-\phi_x\|_{\Cc^0}\leq e^{-9R}$ on $\D_R$ and $\|\psi'-\phi_x\|_{\Cc^2}\leq e^{-6R}$ on $\D_R\setminus\phi_x^{-1}(\cup {1\over 4}\U_e)$. 
Here, we use that $\|\phi_x\|_{\Cc^2}\lesssim e^{2R}$ on $\D_R$ when we consider the Euclidean metric on $\D_R$ and the Hermitian metric on $L_x$.
So, a priori, on $\phi_x^{-1}(\cup {1\over 4}\U_e)$, the Beltrami coefficient associated to this map does not satisfy the condition required for conformally close $(R,e^{-3R})$-points.

Using the maps $\widetilde\Psi_{x,y}$ as in Lemma \ref{lemma_proj_hol}, we can correct $\psi'$ in each connected component of $\phi_x^{-1}({3\over 4}\U_e)$ in order to obtain a map $\psi''$ such that $\|\psi''-\phi_x\|_{\Cc^0}\leq e^{-8R}$ on $\D_R$ and
$\|\psi''-\phi_x\|_{\Cc^2}\leq e^{-5R}$ on $\D_R\setminus\phi_x^{-1}(\cup {1\over 4}\U_e)$. Moreover, $\psi''$ is holomorphic on $\phi_x^{-1}(\cup {1\over 2}\U_e)$. Therefore, its Beltrami's coefficient vanishes on 
$\phi_x^{-1}(\cup {1\over 2}\U_e)$ and hence
satisfies the required property. It remains to modify $\psi''$ in order to obtain a map $\psi$ with $\psi(0)=y$. But this can be done by composition with an automorphism of $\D$, close to the identity, as in the end of the proof of Proposition \ref{prop_conformal_close} or in the proof of Theorem 2.1 in \cite{DinhNguyenSibony2}.
\endproof

The following lemma together with Lemma \ref{lemma_end_1} completes the proof of Proposition \ref{prop_end}.

\begin{lemma} \label{lemma_end_2}
Let $x,y,D$ be as in Proposition \ref{prop_end}. Assume that the tree $F_D$ is not complete. Then, $x$ and $y$ are $(R,e^{-R})$-close.
\end{lemma}
\proof
Since $F_D$ is not complete, there is a path $(\xi^0,\ldots,\xi^m)$ 
of the graph $F_x$ joining $\xi^0=0$ to a vertex $\xi^m$ such that $\phi_x(\xi^m)$ belongs to an exceptional transversal or to $2D'$ with $D'$ an exceptional disc. Define $x^i:=\phi_x(\xi^i)$. 

The image of $(\xi^0,\ldots,\xi^m)$ by $\sigma$ is a path $(\zeta^0,\ldots,\zeta^m)$ of $F_y$. 
We have seen in the proof of Lemma \ref{lemma_discrete_leaf} that $\dist_P(\xi^i,\xi^{i+1})$ and 
$\dist_P(\zeta^i,\zeta^{i+1})$ are smaller than $m_1\hbar\ll \epsilon_0$.
Define $y^i:=\phi_y(\zeta^i)$. There is a disc $D_i$ in $\Vc_{N-i}$ such that $x^i$ and $y^i$ belong to $2D_i$. So, there is a basic projection associated to $x^i$ and $y^i$. Denote by $z^i$ the image of  $x^i$ by this projection.

Observe that $\dist(x^m, z^m)\leq 10\alpha_1$. Recall that $\alpha_1:=e^{-e^{7\lambda R}}$. Therefore, using
Proposition \ref{prop_conformal_close} applied to $\lambda$ large enough, we obtain that $x^m$ and $z^m$ are conformally $(2Nm_1\hbar, e^{-4R})$-close. Moreover, the map $\psi$ associated to these conformally close points is obtained using basic projections as in (H1), (H1)' and Lemma \ref{lemma_proj_E}. So, we can follow the path $(\xi^0,\ldots, \xi^m)$ and see that $\xi^i$ is sent by $\psi$ to $z^i$. So, $0$ is sent by $\psi$ to $z:=z^0$. By Lemma \ref{lemma_discrete_leaf}, the Poincar\'e distance between 0 and $\xi^m$ is at most equal to $mm_1\hbar$. Therefore, $x$ and $z$ are conformally $(Nm_1\hbar,e^{-4R})$-close. It follows that $x$ and $z$ are $(R,e^{-2R})$-close. Moreover, $z$ belongs to a small plaque containing $y$ and $\dist(z,y)\leq 2\dist(x,y)\leq e^{-3R}$. We deduce that $z$ and $y$ are $(R,e^{-2R})$-close. This implies that $x$ and $y$ are $(R,e^{-R})$-close.
\endproof

\noindent
T.-C. Dinh, UPMC Univ Paris 06, UMR 7586, Institut de
Math{\'e}matiques de Jussieu, 4 place Jussieu, F-75005 Paris,
France.\\
{\tt  dinh@math.jussieu.fr}, {\tt http://www.math.jussieu.fr/$\sim$dinh}

\medskip

\noindent
V.-A.  Nguy{\^e}n, 
Vietnamese Academy of Science and Technology, Institute of Mathematics, Department of Analysis, 18 Hoang Quoc Viet Road, Cau Giay District, 10307 Hanoi, Vietnam. {\tt nvanh@math.ac.vn}\\
Current address: Math{\'e}matique-B{\^a}timent 425, UMR 8628, Universit{\'e} Paris-Sud,
91405 Orsay, France.\\
 {\tt VietAnh.Nguyen@math.u-psud.fr}, {\tt http://www.math.u-psud.fr/$\sim$vietanh}

\medskip

\noindent
N. Sibony, Math{\'e}matique-B{\^a}timent 425, UMR 8628, Universit{\'e} Paris-Sud,
91405 Orsay, France.\\
{\tt Nessim.Sibony@math.u-psud.fr}


\small

\begin{thebibliography}{99}

 
 

\bibitem{DinhNguyenSibony2}
Dinh T.-C., Nguyen V.-A.  and Sibony N.,  Entropy for hyperbolic Riemann surface laminations I,
{\tt arXiv:1105.2307}
 
\bibitem{FornaessSibony2}
Forn\ae ss J.-E. and Sibony N.,   Riemann surface laminations with singularities, 
{\it J. Geom. Anal.}, {\bf 18} (2008), no. 2, 400-442.

 
 
\bibitem{Glutsyuk} Glutsyuk A.A.,  Hyperbolicity of the leaves of a generic one-dimensional holomorphic foliation on a nonsingular projective algebraic variety. (Russian)
 {\it Tr. Mat. Inst. Steklova}, {\bf 213} (1997), Differ. Uravn. s Veshchestv. i Kompleks. Vrem., 90-111; translation in {\it Proc. Steklov Inst. Math.} 1996, no. 2, {\bf 213}, 83-103. 
 
 
 
\bibitem{Neto}
Lins Neto A.,  Uniformization and the Poincar\'e metric on the leaves of a foliation by curves,
{\it  Bol. Soc. Brasil. Mat. (N.S.),}  {\bf 31} (2000), no. 3, 351-366. 

 \bibitem{NetoSoares}
Lins Neto A.  and Soares M.G.,
Algebraic solutions of one-dimensional foliations,
{\it J. Differential Geom.,}  {\bf 43} (1996), no. 3, 652-673. 

\bibitem{Shcherbakov}
Shcherbakov A.A., Metrics and smooth uniformisation of leaves of holomorphic foliations, {\it Mosc. Math. J.}, {\bf 11:1} (2011), 157-178.

\end{thebibliography}
\end{document}